\documentclass{article}
\usepackage{latexsym}
\usepackage{epsfig}
\usepackage{amsmath}
\usepackage{amssymb}
\usepackage{amsthm}
\usepackage{graphicx}
\usepackage{caption}
\usepackage{subcaption}
\newtheorem{theorem}{Theorem}[section]
\newtheorem{corollary}[theorem]{Corollary}

\newtheorem{lemma}[theorem]{Lemma}

\theoremstyle{remark}

\theoremstyle{definition}
\newtheorem{definition}[theorem]{Definition}

\newcommand{\R}{{\mathbb R}}

\DeclareMathOperator\arccot{arccot}

\begin{document}

\title{Extremal cubics on the circle and the 2-sphere}

\author{Anastasia Ivanova\thanks{%
Univ.\ Grenoble Alpes, LJK, 38000 Grenoble, France, HSE University, Moscow, Russian Federation
({\tt anastasia.ivanova@univ-grenoble-alpes.fr}).}
\and
Roland Hildebrand \thanks{%
Univ.\ Grenoble Alpes, CNRS, Grenoble INP, LJK, 38000 Grenoble, France
({\tt roland.hildebrand@univ-grenoble-alpes.fr}).}}

\maketitle

\begin{abstract}
We study balls of homogeneous cubics on $\mathbb R^n$, $n = 2,3$, which are bounded by unity on the unit sphere. For $n = 2$ we completely describe the facial structure of this norm ball, while for $n = 3$ we classify all extremal points and describe some families of faces.
\end{abstract}

\section{Introduction}

The problem of maximizing a homogeneous cubic over the unit sphere is NP-hard and arises in various applications in non-convex and combinatorial optimization, e.g., the stable set problem \cite{NesterovMatEllipsoid}. Several approaches for the solution of this problem have been proposed \cite{ManChoSo11,Nie12,ZhangQiYe12,AhmedStill19,BuchheimFampaSarmiento19,deKlerkLaurent20online,FangFawzi21}.

The problem is equivalent to determining the intersection of a given ray in the space of homogeneous cubics with the boundary of the unit norm ball $B_1(S^{n-1})$ in this space, where the norm is defined by
\begin{equation} \label{def_norm}
\|p\| = \max_{\|x\|_2 = 1}\, p(x),
\end{equation}
and $\|\cdot\|_2$ is the usual Euclidean norm in $\mathbb R^n$. A semi-definite description of this norm ball for $n = 3$ follows from Propositions 4.8 and 5.1 in \cite{Saunderson19}, for $n = 2$ such a description is readily constructed from the semi-definite description of the cone of nonnegative univariate trigonometric polynomials \cite[Section 3.4]{NesterovSOS}. Inequivalent semi-definite descriptions which can be generalized to balls of non-homogeneous cubics can be found in \cite{Hildebrand21b}.

In this work we analyze the corresponding semi-definite representable norm balls in the space of homogeneous cubics, which have dimensions 10 and 4, respectively. In particular, we determine all their extreme points and study their facial structure.

Extremal nonnegative polynomials have among others been studied in \cite{ChoiLam77,Reznick78,Reznick00,Naldi14}. The facial structure of cones of nonnegative polynomials has been studied in \cite{KunertThesis,BlekhermanEtAl15}. In particular, the exposed faces of the cone of nonnegative polynomials have been characterized as the set of nonnegative polynomials which vanish on a certain subset of points \cite[Section 4.4.5]{BlekhermanParriloBook}, \cite[Prop.~1.69]{KunertThesis}.

Our motivation for studying the structure of the norms balls $B_1(S^{n-1})$ for $n = 2,3$ is that in the analysis of optimization algorithms on self-concordant functions by methods of optimal control these balls appear as control sets, and the optimal controls are extreme points of these balls \cite{Hildebrand21}.

The remainder of the paper is structured as follows. In Section \ref{sec:definitions} we define some basic notations, in particular, faces of a convex set and investigate the connection between the faces of the norm ball $B_1(S^{n-1})$ and the maxima of cubics on $S^{n-1}$. In Section \ref{sec:n2} we consider the norm ball $B_1(S^1)$, while in Section \ref{sec:n3} we consider $B_1(S^2)$.

\section{Definitions and preliminary results} \label{sec:definitions}

We consider spaces of homogeneous cubic polynomials in $\mathbb R^n$, $n = 2,3$. We call a point $x \in S^{n-1}$ a \emph{global maximum} of the cubic $p$ if $p(x) = \|p\|$. In particular, for any $p \in \partial B_1(S^{n-1})$ the polynomial assumes the value 1 at its global maxima.

The main object of study is the unit norm ball $B_1(S^{n-1})$ in these spaces, i.e., the set of homogeneous cubics which are bounded by 1 on the unit sphere $S^{n-1}$ (equivalently, on the unit ball in $\mathbb R^n$). If $U: \mathbb R^n \to \mathbb R^n$ is an orthogonal transformation, then $B_1(S^{n-1})$ is invariant under the linear map $p \mapsto q$, where $q(x) = p(U^{-1}x)$. This definition is motivated as follows. If $x \in S^{n-1}$ is a global maximum of $p$, then $Ux$ is a global maximum of $q$ on $S^{n-1}$. For ease of notation we shall say that $U$ maps $p$ to $q$, or $q$ is the image of $p$ under $U$.

We call a cubic polynomial $p$ \emph{zonal} if there exists a linear functional $l$ on $\mathbb R^n$ and real constants $\lambda_1,\lambda_3$ such that $p(x) = \lambda_1 \cdot l(x)\cdot \|x\|_2^2 + \lambda_3 \cdot l(x)^3$. In other words, the values of $p$ on the unit sphere $S^{n-1}$ depend only on the values of $l$.

An important role will be played by the faces of $B_1(S^{n-1})$.

{\definition Let $X \subset V$ be a convex set in some real vector space $V$. A \emph{face} of $X$ is a subset $F \subset X$ such that if $l \subset X$ is a line segment and $p \in F$ is an interior point of $l$, then $l \subset F$. A point $p \in X$ is called \emph{extremal} if $\{p\}$ is a face of $X$, or equivalently, if there does not exist a line segment $l \subset X$ such that $p$ is an interior point of $l$. }

The intersection of faces is again a face, and therefore for every $p \in X$ there exists a smallest face $F$ of $X$ which contains $p$. We shall call this face the \emph{minimal face} generated by $p$ and denote it by $F_p$. In particular, $p$ is extremal in $X$ if $\{p\}$ is its minimal face. In this case $p$ is also extremal in every face $F$ containing it.

Our strategy to classify the extreme points of $B_1(S^{n-1})$ is to descend the chain of faces from the larger-dimensional ones to the smaller-dimensional ones up to the zero-dimensional extreme points. Smaller faces are characterized by larger sets of constraints on the polynomials contained in these faces. 

Just as the faces of cones of nonnegative polynomials are defined by constraints imposed by zeros of these polynomials, in our case the faces of the norm ball $B_1(S^{n-1})$ of polynomials are defined by constraints imposed by maxima of these polynomials on the sphere $S^{n-1}$. 
The following result makes the imposed constraints more explicit.

{\lemma \label{lem:order_face} Let $\gamma: \mathbb R \to S^{n-1}$ be an analytic curve, and let $k \geq 0$ be an integer. For every $p \in B_1(S^{n-1})$, define the nonnegative analytic function $f_p(t) = 1 - p(\gamma(t))$ on $\mathbb R$. Then the set $F = \{ p \in B_1(S^{n-1}) \mid f_p(t) = o(t^k),\ t \to 0 \}$ is a face of $B_1(S^{n-1})$. }

\begin{proof}
Let $p \in F$, and let $l \subset B_1(S^{n-1})$ be a line segment containing $p$ in its interior. Then there exists a non-zero cubic $q$ on $\mathbb R^n$ such that $p_{\pm} = p \pm q \in l$. Consider the analytic function $\delta(t) = q(\gamma(t))$.

Suppose for the sake of contradiction that $\delta(t) = ct^m + o(t^m)$ for some integer $m \leq k$ and some non-zero coefficient $c$. Both functions $f_{p_{\pm}}(t) = f_p(t) \mp \delta(t)$ are nonnegative on $\mathbb R$, because $p_{\pm} \in B_1(S^{n-1})$, but on the other hand $f_{p_{\pm}}(t) = \mp ct^m + o(t^m)$, because $f_p(t) = o(t^m)$. This leads to a contradiction, and $\delta(t)$ must be of order $o(t^k)$.

Therefore $f_{p_{\pm}}(t) = o(t^k)$ and hence $p_{\pm} \in F$. It follows that if $F$ contains an interior point of $l$, then it contains also a neighbourhood of that point in $l$. But $F$ is obviously a closed set, and thus contains the whole interval $l$.
\end{proof}

{\corollary \label{cor:minimal_face} Let $p \in B_1(S^{n-1})$, and let $F_p \subset B_1(S^{n-1})$ be the minimal face generated by $p$. Let $\gamma: \mathbb R \to S^{n-1}$ be an analytic curve, and let $k \geq 0$ be an integer. Suppose that the function $f_p(t) = 1 - p(\gamma(t))$ is of order $o(t^k)$ at $t = 0$. Then for every cubic $\tilde p \in F_p$, the function $f_{\tilde p}(t) = 1 - \tilde p(\gamma(t))$ is also of order $o(t^k)$ at $t = 0$. }

\begin{proof}
Since $F_p$ is the minimal face generated by $p$, it is contained in the face $\{ \tilde p \in B_1(S^{n-1}) \mid f_{\tilde p}(t) = o(t^k) \}$ constructed in the previous lemma.
\end{proof}

It follows that for a cubic $p \in \partial B_1(S^{n-1})$, every global maximum $x \in S^{n-1}$ of $p$ imposes a number of linear conditions on the polynomials in the minimal face $F_p$ of $B_1(S^{n-1})$ generated by $p$. Namely, let $\gamma: \mathbb R \to S^{n-1}$ be an analytic curve such that $x = \gamma(0)$. Then the function $f_p(t) = 1 - p(\gamma(t))$ is analytic and nonnegative on $\mathbb R$, and $f(0) = 0$. It follows that either $f(t) \equiv 0$ or there exists $c > 0$ and $k > 0$ such that $f_p(t) = ct^{2k} + o(t^{2k})$. For every polynomial $\tilde p$ in $F_p$ we then have that $\frac{d^m}{dt^m}f_{\tilde p}(t)|_{t = 0} = 0$ for every $m < 2k$, yielding a number of linear conditions on the polynomial $\tilde p$.

The more the maximum $x$ of $p$ is degenerated, i.e., the higher the order of the functions $f_p(t) = 1 - p(\gamma(t))$ defined by curves $\gamma$ through $x = \gamma(0)$, the more linear conditions it induces on the polynomials in $F_p$. If there are enough global maxima to determine $p$ uniquely by the induced linear conditions, then $F_p = \{p\}$ and $p$ is extremal in $B_1(S^{n-1})$.

{\definition \label{def:maxima_types}
In the study of extremal points of $B_1(S^2)$ we will encounter the following types of isolated maxima $x \in S^2$ of a cubic $p$:
\begin{itemize}
\item non-degenerate: the Hessian of $p|_{S^2}$ at $x$ is negative definite;
\item double: the Hessian of $p|_{S^2}$ at $x$ has rank 1, and the order of $f_p(t) = 1 - p(\gamma(t))$ for any non-degenerate analytic curve $\gamma: \mathbb R \to S^2$ which is tangent to the kernel vector of the Hessian at $x = \gamma(0)$ equals 4;
\item triple: the Hessian of $p|_{S^2}$ at $x$ has rank 1, and there exists a non-degenerate analytic curve through $x = \gamma(0)$ such that the order of $f_p(t) = 1 - p(\gamma(t))$ at $x$ equals 6;
\item flat: the Hessian of $p|_{S^2}$ at $x$ vanishes.
\end{itemize}
}

We shall denote the simplex $\{ x = (x_1,\dots,x_n)^T \in \mathbb R^n \mid x_i \geq 0,\ x_1 + \dots + x_n = 1 \}$ by $\Delta_n$ and its relative interior by $\Delta_n^o$.

\section{Norm ball of cubic polynomials on $S^1$} \label{sec:n2}

In this section we investigate the structure of the unit norm ball $B_1(S^1)$ in the space of homogeneous cubic polynomials on $\mathbb R^2$. In particular, we classify the faces of $B_1(S^1)$, including its extremal points.

Let us introduce polar coordinates $r,\varphi$ in $\mathbb R^2$. Note that $p$ is bounded on $S^1$ from below by $-1$ if and only if it is bounded from above by 1. Hence a polynomial $p(x) = p_{30}x_1^3 + 3p_{21}x_1^2x_2 + 3p_{12}x_1x_2^2 + p_{03}x_2^3$ is in $B_1(S^1)$ if and only if
\[ f(\varphi) = p_{30}\cos^3\varphi + 3p_{21}\cos^2\varphi\sin\varphi + 3p_{12}\cos\varphi\sin^2\varphi + p_{03}\sin^3\varphi + 1 \geq 0 \qquad \forall\ \varphi \in [-\pi,\pi].
\]

We shall now describe the boundary of the ball $B_1(S^1)$. Clearly $p \in \partial B_1(S^1)$ if and only if there exists $\varphi^* \in [-\pi,\pi]$ such that $f(\varphi^*) = f'(\varphi^*) = 0$. By a rotation of the coordinate system in $\mathbb R^2$ we may achieve that $\varphi^* = \pi$. Equivalently, the polynomial $p \in B_1(S^1)$ satisfies $p(-e_1) = -1$ and by homogeneity $p(e_1) = 1$, where $e_1 = (1,0)^T$ is the first unit basis vector in $\mathbb R^2$. By Lemma \ref{lem:order_face} the set of all $p \in B_1(S^1)$ such that $p(e_1) = 1$ forms a face $F$ of $B_1(S^1)$. It is described in the following lemma.

{\lemma \label{lem:face_S1} Let $F$ be the face of the unit ball $B_1(S^1)$ defined by the set of polynomials $p \in B_1(S^1)$ such that $p(e_1) = 1$. Then
\begin{equation} \label{face_F}
\begin{aligned}
F &= \left\{ p(x) = x_1^3 + 3p_{12}x_1x_2^2 + p_{03}x_2^3 \left| p_{12} \geq -1,\ p_{03}^2 \leq 1 - 2p_{12}^3 - 3p_{12}^2 \right. \right\} \\
&= \left\{ p(x) = x_1^3 + 3p_{12}x_1x_2^2 + p_{03}x_2^3 \left| \begin{pmatrix} 3 - 6p_{12} & 2p_{03} & 2p_{12} - 1 \\ 2p_{03} & 2p_{12}+2 & 0 \\ 2p_{12} - 1 & 0 & 1 \end{pmatrix} \succeq 0 \right. \right\}.
\end{aligned}
\end{equation} }

\begin{figure}
\centering
\includegraphics[width=9.98cm,height=8.02cm]{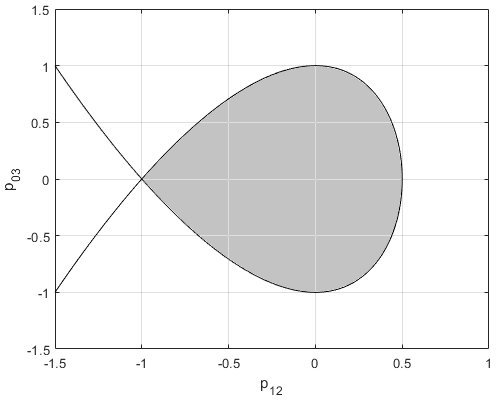}
\caption{Solution curve of the equation $1 - 2p_{12}^3 - 3p_{12}^2 - p_{03}^2 = 0$ delimiting the face $F$.}
\label{fig:faceF}
\end{figure}

\begin{proof}
The conditions $f(\pi) = f'(\pi) = 0$ are equivalent to $p_{30} = 1$, $p_{21} = 0$. We obtain that $p \in F$ if and only if
\[ f(\varphi) = \cos^3\varphi + 3p_{12}\cos\varphi\sin^2\varphi + p_{03}\sin^3\varphi + 1 \geq 0 \qquad \forall\ \varphi \in (-\pi,\pi).
\]
Let us change the independent variable $\varphi$, setting $\cos\varphi = \frac{1 - t^2}{1 + t^2}$, $\sin\varphi = \frac{2t}{1 + t^2}$. The above condition becomes
\begin{equation} \label{g_nonneg_condition}
\begin{aligned}
f(\varphi(t)) = g(t) &= (1 - t^2)^3 + 12p_{12}t^2(1 - t^2) + 8p_{03}t^3 + (1 + t^2)^3 \\ &= 2 \left( (3 - 6p_{12})t^4 + 4p_{03}t^3 + 6p_{12}t^2 + 1 \right) \geq 0 \qquad \forall\ t \in \mathbb R.
\end{aligned}
\end{equation}

Every nonnegative univariate polynomial is a sum of squares of polynomials of lower degree \cite{Hilbert}. This allows to obtain a semi-definite description of the set of nonnegative univariate polynomials \cite[Section 3.1]{NesterovSOS}. In particular, a polynomial $c_{2n}t^{2n} + c_{2n-1}t^{2n-1} + \dots + c_0$ is nonnegative on $\mathbb R$ if and only if there exists a positive semi-definite matrix $C$ of size $n+1$ whose skew-diagonals sum to the coefficients $c_i$ of the polynomial. Thus condition \eqref{g_nonneg_condition} is equivalent to the linear matrix inequality
\begin{equation} \label{F_LMI}
\exists\ \alpha \in \mathbb R: \qquad M(\alpha) = \begin{pmatrix} 3 - 6p_{12} & 2p_{03} & \alpha \\ 2p_{03} & 6p_{12}-2\alpha & 0 \\ \alpha & 0 & 1 \end{pmatrix} \succeq 0.
\end{equation}

The matrix $M(\alpha)$ defines a 3-dimensional spectrahedron in the space of variables $p_{12},p_{03},\alpha$, of which the face $F$ is a projection. We shall now show that $F$ can actually be represented as an affine section of this spectrahedron, by specifying $\alpha$ as an affine function of the remaining two variables.

Suppose the pair $(p_{12},p_{03})$ satisfies condition \eqref{F_LMI}. Then $3 - 6p_{12} \geq 0$, and $p_{12} \leq \frac12$. Consider the determinant
\[ D(\alpha) = \det\,M(\alpha) = 2\left( \alpha^3 - 3p_{12}\alpha^2 + (6p_{12} - 3)\alpha - 2p_{03}^2 - 18p_{12}^2 + 9p_{12} \right).
\]
We have $\frac{dD(\alpha)}{d\alpha} = 6(\alpha - 1)(\alpha - 2p_{12} + 1)$, $\frac{d^2D(\alpha)}{d\alpha^2} = 12(\alpha - p_{12})$, and hence $D(\alpha)$ has a strict local maximum at $\alpha = 2p_{12} - 1 \leq 0$ and a strict local minimum at $\alpha = 1$.

The set of values $\alpha \in \mathbb R$ such that $M(\alpha) \succeq 0$ forms a compact interval, at both ends of which the determinant $D(\alpha)$ vanishes and on the interior of which it is positive. Therefore the local maximum $\alpha = 2p_{12} - 1$ is contained in this interval and gives rise to a positive semi-definite matrix $M(\alpha)$. On the other hand, if for arbitrary values of $p_{12},p_{03}$ the matrix $M(\alpha)$ is positive semi-definite at $\alpha = 2p_{12} - 1$, then condition \eqref{F_LMI} is satisfied for these values.

We therefore can replace \eqref{F_LMI} by the equivalent condition $M(2p_{12} - 1) \succeq 0$, which yields the representation in the second line of \eqref{face_F}.


The face $F$ is delimited by the solution curve of the equation $D(2p_{12} - 1) = 4(1 - 2p_{12}^3 - 3p_{12}^2 - p_{03}^2) = 0$. Since $F$ is convex and compact, it must correspond to the shaded region in Fig.~\ref{fig:faceF}. This yields the representation in the first line of \eqref{face_F}.
\end{proof}

{\corollary \label{cor:bdB1} The boundary of the unit norm ball $B_1(S^1)$ in the space of homogeneous cubic polynomials on $\mathbb R^2$ is given by the union $\bigcup_{\varphi \in [-\pi,\pi]} F_{\varphi}$, where $F_{\varphi}$ is the image of face \eqref{face_F} under rotation of the coordinate system in $\mathbb R^2$ by an angle of $\varphi$. }

\begin{proof}
Let $p \in \partial B_1(S^1)$. Then there exists $\varphi \in [-\pi,\pi]$ such that $p(x^*) = 1$ with $x^* = (\cos\varphi,\sin\varphi)^T$. Rotation of the coordinate system by $-\varphi$ maps $x^*$ to the unit basis vector $e_1$, and hence $p$ to a point in the face $F$. It follows that $p \in F_{\varphi}$.

On the other hand, for every $\varphi \in [-\pi,\pi]$ the face $F_{\varphi}$ is contained in the boundary $\partial B_1(S^1)$. This completes the proof.
\end{proof}

Let us now investigate the extreme points of the norm ball $B_1(S^1)$. From Corollary \ref{cor:bdB1} it follows that every extreme point of $B_1(S^1)$ can be obtained from some extreme point of face \eqref{face_F} by a coordinate rotation. On the other hand, the (relative) boundary of face \eqref{face_F} does not contain line segments, and every boundary point of this face is extreme.

The boundary of face \eqref{face_F} happens to be defined by a rational curve, i.e., it can be polynomially parameterized by a real variable $\tau$. In the case at hand it is checked by direct calculation that the parameter can be chosen in the interval $\tau \in \left[ -\frac{\sqrt{3}}{2},\frac{\sqrt{3}}{2} \right]$, with the corresponding polynomial on the boundary of face \eqref{face_F} given by
\begin{equation} \label{extremal_poly}
p(x) = x_1^3 + \frac32(1 - 4\tau^2)x_1x_2^2 + (3\tau - 4\tau^3)x_2^3 = -\frac12(x_1 + 2\tau x_2)(4\tau^2x_2^2 + 4\tau x_1x_2 - 2x_1^2 - 3x_2^2).
\end{equation}
Here the parameter value $\tau = 0$ corresponds to the right-most point of the curve in Fig.~\ref{fig:faceF}, while the values $\tau = \pm\frac{\sqrt{3}}{2}$ correspond to the point of self-intersection.

For $\tau = 0$ and $x = (\cos\varphi,\sin\varphi)^T \in S^1$ polynomial \eqref{extremal_poly} simplifies to
\[ p(x) = x_1^3 + \frac32x_1x_2^2 = -\frac12\cos^3\varphi + \frac32\cos\varphi = 1 - \frac38\varphi^4 + O(\varphi^6),
\]
and the maximum at $x = e_1$ on $S^1$ is degenerated (has a vanishing second derivative with respect to $\varphi$).

For $\tau = \pm\frac{\sqrt{3}}{2}$ we get
\[ p(x) = x_1^3 - 3x_1x_2^2 = 4\cos^3\varphi - 3\cos\varphi = \cos3\varphi,
\]
and $p$ has three equally spaced global maxima on $S^1$.

It is easily checked that polynomial \eqref{extremal_poly} is invariant with respect to the reflection
\[ x \mapsto \frac{1}{1 + 4\tau^2} \begin{pmatrix} 1 - 4\tau^2 & 4\tau \\ 4\tau & -1 + 4\tau^2 \end{pmatrix} x.
\]
Hence for $\tau \in \left( -\frac{\sqrt{3}}{2},\frac{\sqrt{3}}{2} \right) \setminus \{0\}$ it has two global maxima on $S^1$, which are an angle $\phi = \arccos\frac{1-4\tau^2}{1+4\tau^2} \in \left(0,\frac{2\pi}{3} \right)$ apart.

We obtain the following result.

{\theorem \label{thm:structure_S1} Let $p \in B_1(S^1)$ be a homogeneous cubic on $\mathbb R^2$ of norm not exceeding 1. Then exactly one of the following conditions holds:
\begin{itemize}
\item $p$ is strictly smaller than 1 on the unit circle $S^1$, then it is contained in the interior of $B_1(S^1)$;
\item $p$ has a single non-degenerate global maximum on $S^1$ with value 1, then it lies in the (relative) interior of the face $F_{\varphi}$ for some $\varphi$;
\item $p$ has a single degenerate global maximum on $S^1$ with value 1, then it is the image of polynomial \eqref{extremal_poly} with parameter value $\tau = 0$ under a rotation of the coordinate system;
\item $p$ has exactly two non-degenerate global maxima on $S^1$ with value 1, then it is the image of polynomial \eqref{extremal_poly} with parameter value $\tau \in \left( -\frac{\sqrt{3}}{2},\frac{\sqrt{3}}{2} \right) \setminus \{0\}$ under a rotation of the coordinate system, and the two maxima are an angle $\phi = \arccos\frac{1-4\tau^2}{1+4\tau^2} \in (0,\frac{2\pi}{3})$ apart;
\item $p$ has exactly three non-degenerate global maxima on $S^1$ with value 1, then it is the image of polynomial \eqref{extremal_poly} with parameter value $\tau = \frac{\sqrt{3}}{2}$ under a rotation of the coordinate system, and the three maxima are equally spaced.
\end{itemize}
In the last three cases the polynomial $p$ is extremal in $B_1(S^1)$.
}

\begin{proof}
The first item follows from the definition of the norm.

Suppose $p$ has a single non-degenerate global maximum with value 1. Then $p \in \partial B_1(S^1)$, and by Corollary \ref{cor:bdB1} it belongs to some face $F_{\varphi}$. Without loss of generality, we may assume $p \in F$, where $F$ is given by \eqref{face_F}. But then $p$ cannot belong to the boundary of $F$, because the boundary points of $F$ either have several global maxima or the unique global maximum is degenerated. Hence $p$ is in the (relative) interior of $F$.

On the other hand, suppose that $p$ is in the interior of some face $F_{\varphi}$, without loss of generality in the interior of $F$. For the sake of contradiction, suppose $p$ has multiple global maxima, or a degenerate global maximum. By Corollary \ref{cor:minimal_face} it shares this property with all elements of $F$, in particular, with $\hat p(x) = x_1^3 + x_1x_2^2$. This polynomial reduces to $\hat p(x) = \cos\varphi$ on $x = (\cos\varphi,\sin\varphi)^T \in S^1$, which has only one non-degenerate global maximum. This leads to a contradiction, proving the second item.

If $p$ does not satisfy one of the first two conditions, then it has to lie on the boundary of some face $F_{\varphi}$. Hence the last three items correspond to the images of polynomial \eqref{extremal_poly} under rotations of the coordinate system, and have been classified above.
\end{proof}

\section{Norm ball of cubic polynomials on $S^2$} \label{sec:n3}

In this section we describe the extremal elements of the unit norm ball $B_1(S^2)$ of homogeneous cubic polynomials on $\mathbb R^3$, where the norm is given by \eqref{def_norm}. Since the norm ball is invariant with respect to orthogonal transformations of the space $\mathbb R^3$, we shall use such maps to transform the polynomials into a convenient form.

Clearly any extremal element $p$ of $B_1(S^2)$ has to lie on the boundary of $B_1(S^2)$, and hence there exists a global maximum $x^* \in S^2$ of $p$ on the unit sphere with value 1. By a rotation of $\mathbb R^3$ we may achieve that $p(x^*) = e_1 = (1,0,0)^T$. Since the gradient of $p$ at $x = e_1$ is proportional to $e_1$, $p$ cannot contain terms which are quadratic in $x_1$ and has to be of the form
\begin{equation} \label{e1face}
p(x) = x_1^3 + 3\left(p_{120}x_2^2 + 2p_{111}x_2x_3 + p_{102}x_3^2\right)x_1 + \left(p_{030}x_2^3 + 3p_{021}x_2^2x_3 + 3p_{012}x_2x_3^2 + p_{003}x_3^3\right).
\end{equation}
Define the functions
\begin{equation} \label{p12_p03_functions}
\begin{aligned}
p_{12}(\varphi) &= p_{120}\cos^2\varphi + 2p_{111}\sin\varphi\cos\varphi + p_{102}\sin^2\varphi, \\
p_{03}(\varphi) &= p_{030}\cos^3\varphi + 3p_{021}\cos^2\varphi\sin\varphi + 3p_{012}\cos\varphi\sin^2\varphi + p_{003}\sin^3\varphi.
\end{aligned}
\end{equation}

{\lemma \label{lem:face_S2} Assume above notations. A polynomial $p(x)$ of the form \eqref{e1face} is in the unit norm ball $B_1(S^2)$ if and only if
\[ p_{12}(\varphi) \geq -1,\ 1 - 2p_{12}(\varphi)^3 - 3p_{12}(\varphi)^2 - p_{03}(\varphi)^2 \geq 0\quad \forall\ \varphi \in [-\pi,\pi]
\]
or, equivalently,
\[ \begin{pmatrix} 3 - 6p_{12}(\varphi) & 2p_{03}(\varphi) & 2p_{12}(\varphi) - 1 \\ 2p_{03}(\varphi) & 2p_{12}(\varphi)+2 & 0 \\ 2p_{12}(\varphi) - 1 & 0 & 1 \end{pmatrix} \succeq 0\quad \forall\ \varphi \in [-\pi,\pi].
\]
The set ${\cal F}$ of all such polynomials forms a 7-dimensional face of $B_1(S^2)$. }

\begin{proof}
Introduce polar coordinates $r,\varphi$ such that $x_2 = r\cos\varphi$, $x_3 = r\sin\varphi$. Then \eqref{e1face} takes the form
\[ p(x) = x_1^3 + 3p_{12}(\varphi)x_1r^2 + p_{03}(\varphi)r^3.
\]
For every fixed $\varphi$ this is a homogeneous cubic polynomial in the two variables $x_1,r$, denote this polynomial by $\tilde p_{\varphi}$.

Clearly $p \in B_1(S^2)$ if and only if $\tilde p_{\varphi} \in B_1(S^1)$ for all $\varphi$. The first assertion of the lemma is now a direct consequence of Lemma \ref{lem:face_S1}.

Let us show the second assertion. The set ${\cal F}$ consists of all polynomials $p \in B_1(S^2)$ such that $p(e_1) = 1$, and hence is even an exposed face of $B_1(S^2)$, namely the intersection of $B_1(S^2)$ with the supporting hyperplane defined by the linear functional $p \mapsto p(e_1)$. Further, polynomial \eqref{e1face} has 7 coefficients whose values are not fixed, and ${\cal F}$ is at most 7-dimensional.

Consider the trigonometric polynomial
\begin{equation} \label{Delta_polynom}
\Delta(\varphi) = 1 - 2p_{12}(\varphi)^3 - 3p_{12}(\varphi)^2 - p_{03}(\varphi)^2.
\end{equation}
For $\hat p(x) = x_1(x_1^2+x_2^2+x_3^2)$, which is of the form \eqref{e1face}, we get $p_{12}(\varphi) \equiv \frac13$, $\Delta(\varphi) \equiv \frac{16}{27}$. Hence the inequalities in the formulation of the lemma are all strict, and remain valid if the seven free coefficients are changed by a small enough disturbance. It follows that ${\cal F}$ contains a 7-dimensional ball around $\hat p$, and its dimension equals 7.
\end{proof}

If for some cubic $p \in {\cal F}$ we have $\Delta(\varphi) > 0$ for all $\varphi$, then also $p_{12}(\varphi) > -1$ for all $\varphi$, and $p$ is contained in the interior of ${\cal F}$. Hence if $p \in {\cal F}$ is extremal, then $\Delta(\varphi)$ must have roots, and the closed curve $\{ (p_{12}(\varphi),p_{03}(\varphi)) \mid \varphi \in [-\pi,\pi] \}$ defined by functions \eqref{p12_p03_functions} must touch the boundary of the shaded region in Fig.~\ref{fig:faceF}. In this case, by a rotation of $\mathbb R^3$ around the $e_1$-axis we may achieve that $\Delta(0) = 0$. Then the restriction of $p$ to the $(e_1,e_2)$-plane is given by \eqref{extremal_poly}.

This condition singles out a sub-face $F \subset {\cal F}$, but its shape depends on the location of the point $(p_{12}(0),p_{03}(0))$ on the curve in Fig.~\ref{fig:faceF}. The different cases which arise are studied in Section \ref{sec:three_maxima} and subsequent sections.

\subsection{Set of global maximizers on $S^2$}

In this section we derive some constraints on the set of global maxima on $S^2$ of a non-zero homogeneous cubic polynomial $p$ on $\mathbb R^3$. Without loss of generality we may normalize the cubic such that $p \in \partial B_1(S^2)$, and the global maxima take the value 1.

{\lemma \label{lem:max_angle}  Let $p \in \partial B_1(S^2)$. Then the angle between two different global maxima of $p$ on $S^2$ is at most $\frac{2\pi}{3}$. If the angle between two global maxima $u,v$ equals $\frac{2\pi}{3}$, then there are exactly three global maxima on the great circle through $u,v$. }

\begin{proof}
Let $H \subset \mathbb R^3$ be the plane spanned by two different global maxima of $p$ on $S^2$. Then the restriction of $p$ to $H$ has two different global maxima on $S^1$, and satisfies either condition 4 or condition 5 in Theorem \ref{thm:structure_S1}. In neither case the maxima can be farther apart than by an angle $\frac{2\pi}{3}$, and equality is attained when condition 5 holds.
\end{proof}

{\lemma \label{lem:exactly3maxima} Let $p \in \partial B_1(S^2)$. If $p$ has at least three distinct global maxima on some great circle, then it has exactly three maxima on this great circle, and they are equally spaced. }

\begin{proof}
Apply Theorem \ref{thm:structure_S1} to the restriction of $p$ to the plane containing the great circle in question.
\end{proof}

{\lemma \label{lem:4_planar_maxima} If $p \in \partial B_1(S^2)$ has at least four distinct global maxima on some circle $C \subset S^2$, then $p$ is a zonal function and constant on $C$. Moreover, it is an extremal point of $B_1(S^2)$. }

\begin{proof}
By the preceding lemma $C$ cannot be a great circle. Without loss of generality, let the centre of $C$ be a multiple of the basis vector $e_1$, i.e., $C = \{ x \in S^2 \mid x_1 = c \}$ for some constant $c \in (0,1)$. Parameterize $C$ by an angle $\varphi \in [-\pi,\pi]$. Then the restriction of $p$ on $C$ is a trigonometric polynomial $s(\varphi)$ of degree at most three.

Let us show that $p$ is actually constant on $C$. Suppose for the sake of contradiction that $s(\varphi) \not\equiv const$. Then the global maxima of $s(\varphi)$ are isolated. Since by assumption there are at least 4 such maxima, the derivative $\frac{ds(\varphi)}{d\varphi}$ changes sign at least 8 times. However, this derivative is also a trigonometric polynomial of degree at most three, leading to a contradiction.

It follows that $s \equiv 1$, and the whole circle $C$ consists of global maxima of $p$. But then on every great circle containing $e_1$ there are at least two global maxima of $p$, which are apart by a constant angle of $2\arccos\,c$. By Theorem \ref{thm:structure_S1} this angle completely determines the restriction of $p$ on the plane containing the great circle in question. In particular, the polynomial $p$ is invariant under rotations about the axis $e_1$, and hence zonal. Since it is completely determined by the conditions established as in Section \ref{sec:definitions} by its set of global maxima, its minimal face coincides with $\{p\}$ and $p$ must be extremal in $B_1(S^2)$.
\end{proof}

{\lemma \label{lem:min4maxima} Let $p \in \partial B_1(S^2)$ be a cubic with only isolated, non-degenerate global maxima. If $p$ has no more than three maxima, then it is not an extreme point of $B_1(S^2)$. }

\begin{proof}
Let $a,b,c \in S^2$ be distinct points such that all global maxima of $p$ are among them. Consider the polynomial
\[ \delta(x) = \langle a \times b,x \rangle \cdot \langle b \times c,x \rangle \cdot \langle c \times a,x \rangle,
\]
where $\times$ denotes the cross product and $\langle \cdot,\cdot \rangle$ the usual scalar product in $\mathbb R^3$. Note that $\delta$ is a non-zero homogeneous cubic on $\mathbb R^3$. Moreover, its value and gradient vanish at $x = a,b,c$.

We now show that there exists $\epsilon > 0$ such that $p \pm \epsilon \delta \in B_1(S^2)$. Indeed, since the Hessian of the restriction $p|_{S^2}$ is negative definite at $x = a,b,c$, the polynomials $p \pm \epsilon\delta$ still have isolated global maxima at $x = a,b,c$ if $\epsilon$ is small enough. Hence there exist neighbourhoods $U_a,U_b,U_c \subset S^2$ of $a,b,c$, respectively, such that $p \pm \epsilon\delta$ do not exceed 1 on $U = U_a \cup U_b \cup U_c$. On the other hand, on $S^2 \setminus U$ the polynomial $\delta$ is bounded, while $p$ does not exceed $1 - \varepsilon$ for some $\varepsilon > 0$. Hence for small enough $\epsilon$ the polynomials $p \pm \epsilon\delta$ also do not exceed 1 on $S^2 \setminus U$.

It follows that $p \pm \epsilon\delta \in B_1(S^2)$. In particular, every face of $B_1(S^2)$ which contains $p$ must also contain the line segment between $p \pm \epsilon\delta$, and $p$ cannot be extremal.
\end{proof}

Next we describe the extremal points of $B_1(S^2)$ corresponding to some special cases.

\subsection{Three equally spaced maxima on a great circle} \label{sec:three_maxima}

In this section we consider the case when the restriction of $p$ to some 2-dimensional subspace of $\mathbb R^3$ satisfies condition 5 in Theorem \ref{thm:structure_S1}.

Return to the notations at the beginning of Section \ref{sec:n3}. Suppose that for $p \in {\cal F}$ we have that $p_{12}(0) = -1$. Then the restriction of $p$ to the plane spanned by the basis vectors $e_1,e_2$ satisfies the fifth condition in Theorem \ref{thm:structure_S1}. From \eqref{extremal_poly} it follows with $\tau = \frac{\sqrt{3}}{2}$ that  $p_{030} = 0$, $p_{120} = -1$, and $p$ has additional global maxima at $\left(\cos\frac{2\pi}{3},\pm\sin\frac{2\pi}{3},0\right)^T$ on $S^2$. At $\varphi = 0$ the curve $(p_{12}(\varphi),p_{03}(\varphi))$ defined by \eqref{p12_p03_functions} passes through the corner $(-1,0)$ of the shaded region in Fig.~\ref{fig:faceF}. Since it lies entirely in this region, the derivatives $\frac{dp_{12}}{d\varphi}$, $\frac{dp_{03}}{d\varphi}$ vanish at $\varphi = 0$, which yields the first order conditions $p_{111} = p_{021} = 0$. The following result describes all $p \in {\cal F}$ satisfying these conditions.

{\lemma \label{lem:3maxima} A cubic
\begin{equation} \label{p_3maxima}
p(x) = x_1^3 + 3\left(-x_2^2 + p_{102}x_3^2\right)x_1 + \left(3p_{012}x_2 + p_{003}x_3\right)x_3^2
\end{equation}
is an element of ${\cal F}$ if and only if
\begin{equation} \label{triangle_3maxima}
1 + p_{102} \pm \sqrt{3}p_{012} \geq 0,\quad 1 - 2p_{102} \geq 0,
\end{equation}
and
\begin{equation} \label{determinant_3maxima}
(1 + p_{102} + \sqrt{3}p_{012})(1 + p_{102} - \sqrt{3}p_{012})(1 - 2p_{102}) \geq p_{003}^2.
\end{equation}
The set ${\cal F}_3$ of these cubics forms a 3-dimensional face of $B_1(S^2)$. Its extremal elements are exactly those which satisfy \eqref{determinant_3maxima} with equality and satisfy either all or only one inequality in \eqref{triangle_3maxima} strictly. }

\begin{proof}
Let $p$ be of the form \eqref{p_3maxima}. By Lemma \ref{lem:face_S2} we have $p \in {\cal F}$ if and only if
\[ p_{102}\sin^2\varphi \geq - 1 + \cos^2\varphi = -\sin^2\varphi,\quad \Delta(\varphi) \geq 0 \qquad \forall\ \varphi \in [-\pi,\pi].
\]
The first condition is equivalent to the inequality $p_{102} \geq -1$. Polynomial \eqref{Delta_polynom} becomes
\begin{align} \label{Delta_3maxima}
\Delta(\varphi) &= \sin^4\varphi\left( 3(p_{102}^2 + 2p_{102} - 3p_{012}^2 + 1)\cos^2\varphi - 6p_{012}p_{003}\cos\varphi\sin\varphi \right. \nonumber\\
& \left. + (- 2p_{102}^3 - 3p_{102}^2 - p_{003}^2 + 1)\sin^2\varphi \right).
\end{align}
Nonnegativity of this polynomial is hence equivalent to the matrix inequality
\[ \begin{pmatrix} 3((1 + p_{102})^2 - 3p_{012}^2) & -3p_{012}p_{003} \\ -3p_{012}p_{003} & (1 + p_{102})^2(1 - 2p_{102}) - p_{003}^2 \end{pmatrix} \succeq 0,
\]
which in turn is equivalent to nonnegativity of the diagonal elements and the determinant $3(1 + p_{102})^2((1 - 2p_{102})((1 + p_{102})^2 - 3p_{012}^2) - p_{003}^2)$.

Together with the inequality $p_{102} \geq -1$ nonnegativity of the diagonal elements is equivalent to \eqref{triangle_3maxima}, while nonnegativity of the determinant is equivalent to \eqref{determinant_3maxima}. This proves the first assertion of the lemma.

The set ${\cal F}_3$ consists of exactly those cubics $p \in \partial B_1(S^2)$ which have global maxima at $x = e_1$, $x = \left(\cos\frac{2\pi}{3},\pm\sin\frac{2\pi}{3},0\right)^T$ on $S^2$. Hence ${\cal F}_3$ is a face of $B_1(S^2)$, and clearly it is 3-dimensional.

The projection of ${\cal F}_3$ onto the $(p_{102},p_{012})$-plane is an equilateral triangle $T$ with corners located at $(-1,0),(\frac12,\pm\frac{\sqrt{3}}{2})$. Points of ${\cal F}_3$ projecting to the boundary of $T$ satisfy $p_{003} = 0$, while the pre-images of interior points are segments of positive length. At the end-points of these segments the boundary of ${\cal F}_3$ is strictly convex. Hence the extreme points of ${\cal F}_3$ are the pre-images of the corners of $T$ and the points of the boundary surface where $p_{003} \not= 0$. This proves the remaining assertions.
\end{proof}

Let us now consider the extremal points of ${\cal F}_3$ in more detail.

The pre-image of the corner $(p_{102},p_{012}) = (-1,0)$ of triangle $T$ is the cubic $x_1^3 - 3x_1(x_2^2+x_3^2)$. Its values on $S^2$ depend only on $x_1$, and hence it is a zonal polynomial. It has global maxima on $S^2$ at $e_1$ and on a circle given by $x_1 = \cos\frac{2\pi}{3} = -\frac12$.

The pre-images of the other corners are obtained by a rotation of $\mathbb R^3$ by $\pm\frac{2\pi}{3}$ in the plane spanned by $e_1,e_2$.

Let now $p$ be one of the remaining extremal points, satisfying \eqref{triangle_3maxima} strictly and \eqref{determinant_3maxima} with equality. Let us further assume without loss of generality that
\begin{equation} \label{p003_3maxima}
p_{003} = \sqrt{((1 + p_{102})^2 - 3p_{012}^2)(1 - 2p_{102})} > 0,
\end{equation}
which can be achieved by the orthogonal transformation $x_3 \mapsto -x_3$ of $\mathbb R^3$. Then the coefficient at $\cos^2\varphi$ in \eqref{Delta_3maxima} is positive, and $\Delta(\varphi)$ has a unique pair of roots $\varphi^*,\varphi^* - \pi$ different from $0,\pm\pi$. More precisely, we have
\[ \varphi^* = \arccot\frac{p_{012}p_{003}}{(1 + p_{102})^2 - 3p_{012}^2} = \arccot\frac{p_{012}\sqrt{1 - 2p_{102}}}{\sqrt{(1 + p_{102})^2 - 3p_{012}^2}},
\]
\[ \cos\varphi^* = \frac{p_{012}\sqrt{1 - 2p_{102}}}{\sqrt{(1 + p_{102})(1 + p_{102} - 2p_{012}^2)}},\quad \sin\varphi^* = \frac{\sqrt{(1 + p_{102})^2 - 3p_{012}^2}}{\sqrt{(1 + p_{102})(1 + p_{102} - 2p_{012}^2)}}.
\]
Additional global maxima of $p$ can hence only be located on the great circle through $e_1$ corresponding to this angle $\varphi^*$.

Let us compute the restriction of $p$ to this great circle. Set $x_2 = r\cos\varphi^*$, $x_3 = r\sin\varphi^*$, then \eqref{p_3maxima} becomes
\[ p(x) = x_1^3 + 3\frac{p_{102} + p_{102}^2 - p_{012}^2}{1 + p_{102} - 2p_{012}^2}x_1r^2 + \frac{\sqrt{(1 + p_{102})(1 - 2p_{102})}((1 + p_{102})^2 - 3p_{012}^2)}{(1 + p_{102} - 2p_{012}^2)^{3/2}}r^3,
\]
which has the form \eqref{extremal_poly} with $x_2 = r$ and $\tau = \frac{\sqrt{(1 + p_{102})(1 - 2p_{102})}}{2\sqrt{1 + p_{102} - 2p_{012}^2}} \in \left(0,\frac{\sqrt{3}}{2}\right)$. By Theorem \ref{thm:structure_S1}, item 4, this restriction has a single global non-degenerate maximum besides $e_1$, which is an angle $\arccos\frac{1-4\tau^2}{1+4\tau^2}$ apart from $e_1$. Hence this maximum is at
\begin{equation} \label{p_x_relation_3maxima}
\begin{aligned}
x_1 = \frac{1-4\tau^2}{1+4\tau^2} = \frac{p_{102}+p_{102}^2-p_{012}^2}{1-p_{102}^2-p_{012}^2}, &\quad r = \frac{4\tau}{1 + 4\tau^2} = \frac{\sqrt{(1 + p_{102})(1 - 2p_{102})(1 + p_{102} - 2p_{012}^2)}}{1-p_{102}^2-p_{012}^2}, \\
x_2 = \frac{p_{012}(1 - 2p_{102})}{1-p_{102}^2-p_{012}^2}, &\quad x_3 = \frac{\sqrt{(1 - 2p_{102})((1 + p_{102})^2 - 3p_{012}^2)}}{1-p_{102}^2-p_{012}^2}.
\end{aligned}
\end{equation}
It is checked straightforwardly that this maximum is non-degenerate not only on the great circle, but also on $S^2$.

{\corollary \label{cor:3maxima_extremal} Let $p$ be an extremal element of $B_1(S^2)$ such that there exists a great circle in $S^2$ on which $p$ has three equally spaced global maxima. Then exactly one of the following conditions holds:
\begin{itemize}
\item there exists a rotation of $\mathbb R^3$ which takes $p$ to the zonal cubic $x_1^3 - 3x_1(x_2^2+x_3^2)$, in this case the set of global maxima of $p$ on $S^2$ consists of a single point and a circle whose points are an angle $\frac{2\pi}{3}$ apart from this point;
\item there exists an orthogonal transformation of $\mathbb R^3$ which takes $p$ to the cubic \eqref{p_3maxima} with $p_{003}$ given by \eqref{p003_3maxima} and $p_{102},p_{012}$ satisfying \eqref{triangle_3maxima} strictly, in this case $p$ has four isolated non-degenerate global maxima on $S^2$, three of which lie equally spaced on a great circle, and the fourth of which is an angle strictly less than $\frac{2\pi}{3}$ apart from each of the others.
\end{itemize}
On the other hand, suppose the set of global maxima on $S^2$ of a polynomial $p \in \partial B_1(S^2)$ is as described in either of the two cases above. Then $p$ is an extremal point of $B_1(S^2)$. }

\begin{proof}
Let the cubic $p$ be as in the first part of the corollary. Then $p \in \partial B_1(S^2)$, and the global maxima of $p$ on $S^2$ have value 1. We may rotate $\mathbb R^3$ such that the great circle carrying three global maxima comes to lie in the $(e_1,e_2)$-plane and one of the maxima is moved to $e_1$. Then $p$ takes the form \eqref{p_3maxima}, and the first assertion of the corollary follows from Lemma \ref{lem:3maxima} and the considerations above.

Now let $p \in \partial B_1(S^2)$ be arbitrary. Then again the global maximum of $p$ on $S^2$ equals 1.

If the set of global maxima is as described in the first item of Corollary \ref{cor:3maxima_extremal}, then we may rotate the single point to the basis vector $e_1$, and the circle of maxima comes to be given by the relation $x_1 = -\frac12$. But then on every great circle through $e_1$ the image $\tilde p$ of $p$ under the rotation has three equally spaced global maxima, and its restriction to this great circle is given by \eqref{extremal_poly} with $\tau = \frac{\sqrt{3}}{2}$ by Theorem \ref{thm:structure_S1}. It follows that $\tilde p(x) = x_1^3 - 3x_1(x_2^2+x_3^2)$. But this polynomial is extremal in $B_1(S^2)$, and so is the original polynomial $p$. This proves the first part of the corollary.

Now suppose that $p$ has 4 global maxima, located as in the second item of the corollary. Let us rotate the great circle carrying three maxima into the $(e_1,e_2)$-plane and one of the maxima into $e_1$. Assume by possibly applying a reflection of the $x_3$ coordinate that the fourth maximum $x^*$ of the image $\tilde p$ of $p$ under this transformation has a positive last entry $x^*_3 = \sqrt{1 - (x_1^*)^2 - (x_2^*)^2}$. Since $x^*$ is strictly closer than an angle $\frac{2\pi}{3}$ to each of the other three maxima, we have that $(x_1^*,x_2^*)$ lies in the interior of the triangle $T_x$ spanned by $(1,0)$, $\left(-\frac12,\pm\frac{\sqrt{3}}{2}\right)$. Let $T_p$ be the triangle of pairs $(p_{102},p_{012})$ satisfying \eqref{triangle_3maxima}. We now show that the map
\[ \begin{pmatrix} p_{102} \\ p_{012} \end{pmatrix} \mapsto \begin{pmatrix} x_1 \\ x_2 \end{pmatrix} = \frac{1}{1-p_{102}^2-p_{012}^2}\begin{pmatrix} p_{102}+p_{102}^2-p_{012}^2 \\ p_{012}(1 - 2p_{102}) \end{pmatrix}
\]
derived in \eqref{p_x_relation_3maxima} is a bijection between the interiors of $T_p$ and $T_x$.

Introduce barycentric coordinates $\lambda = (\lambda_1,\lambda_2,\lambda_3)^T$ in $T_p$ and $\mu = (\mu_1,\mu_2,\mu_3)^T$ in $T_x$, such that
\[ \begin{pmatrix} p_{102} \\ p_{012} \end{pmatrix} = \begin{pmatrix} -1 & \frac12 & \frac12 \\ 0 & -\frac{\sqrt{3}}{2} & \frac{\sqrt{3}}{2} \end{pmatrix}\lambda, \quad \begin{pmatrix} x_1 \\ x_2 \end{pmatrix} = \begin{pmatrix} 1 & -\frac12 & -\frac12 \\ 0 & \frac{\sqrt{3}}{2} & -\frac{\sqrt{3}}{2} \end{pmatrix}\mu.
\]
Rewriting the relation between $(p_{102},p_{012})$ and $(x_1,x_2)$ in these coordinates we obtain
\[ \mu = \frac{1}{\lambda_1^{-1} + \lambda_2^{-1} + \lambda_3^{-1}}\begin{pmatrix} \lambda_1^{-1} \\ \lambda_2^{-1} \\ \lambda_3^{-1} \end{pmatrix},
\]
which shows the desired bijection.

Thus there exists an extremal point $\hat p$ of $B_1(S^2)$ of the form \eqref{p_3maxima} with $p_{003}$ given by \eqref{p003_3maxima} which has the same global maxima on $S^2$ as $\tilde p$. However, these maxima determine $\hat p$ uniquely. This implies that $\tilde p$ actually coincides with $\hat p$ and is extremal. But then also $p$ is extremal. This proves the second part of the corollary.
\end{proof}


\subsection{Degenerate maximum}

In this section we consider the case when the restriction of $p$ to some 2-dimensional subspace of $\mathbb R^3$ satisfies condition 3 in Theorem \ref{thm:structure_S1}, or equivalently, if $p$ has a degenerate global maximum on $S^2$.

Suppose that for $p \in {\cal F}$ we have that $p_{12}(0) = \frac12$. From \eqref{extremal_poly} with $\tau = 0$ it follows that $p_{030} = 0$, $p_{120} = \frac12$. Moreover, the function $p_{12}(\varphi)$ defined by \eqref{p12_p03_functions} achieves a global maximum at $\varphi = 0$, which yields the first order condition $\frac{dp_{12}}{d\varphi}|_{\varphi = 0} = 2p_{111} = 0$. The following result describes the corresponding polynomials.

{\lemma \label{lem:degenerated_maxima} The set ${\cal F}_4$ of cubics $p \in B_1(S^2)$ of the form
\begin{equation} \label{p_degenerate_maximum}
p(x) = x_1^3 + 3\left(\frac12x_2^2 + p_{102}x_3^2\right)x_1 + \left(3p_{021}x_2^2x_3 + 3p_{012}x_2x_3^2 + p_{003}x_3^3\right)
\end{equation}
forms a 4-dimensional face of $B_1(S^2)$. Let $p$ be an extremal point of ${\cal F}_4$. Then exactly one of the following conditions holds:
\begin{itemize}
\item $p(x) = x_1^3 + 3\left(\frac12x_2^2 - x_3^2\right)x_1 \pm \frac{3\sqrt{3}}{2}x_2^2x_3$ is zonal, in this case condition 1 of Corollary \ref{cor:3maxima_extremal} holds;
\item $p(x) = x_1^3 + \frac32\left(x_2^2 + 2p_{102}x_3^2\right)x_1 \pm \sqrt{1-2p_{102}}\left(\frac32x_2^2x_3 + (1+p_{102})x_3^3\right)$ for some $p_{102} \in \left(-1,\frac12\right)$ is zonal, in this case the global maximum on $S^2$ is achieved on some circle;
\item $p(x) = x_1^3 + \frac32x_1(x_2^2+x_3^2)$ is zonal, in this case $e_1$ is the unique global maximum of $p$ on $S^2$, and the Hessian of $p|_{S^2}$ identically vanishes at $x = e_1$;
\item $p(x) = x_1^3 + 3\left(\frac12x_2^2 + p_{102}x_3^2\right)x_1 \pm \sqrt{1-2p_{102}}\left(-\frac32x_2^2 + (1+p_{102})x_3^2\right)x_3$ for some $p_{102} \in \left(-1,\frac12\right)$, in this case $p$ has a degenerate global maximum at $e_1$ such that $(1-p)|_{C^{\mp}}$ has order of smallness 6 at $x = e_1$, where $C^{\mp}$ is given by \eqref{circles_triply_degenerated}, and a non-degenerate global maximum at $x^{\pm} = \frac{1}{1 - p_{102}}\left( p_{102}, 0, \pm\sqrt{1-2p_{102}} \right)^T$;
\item $p(x)$ is given by \eqref{deg_nondeg_nondeg} for some $\xi \in [-\pi,\pi] \setminus \{\pm\frac{\pi}{2}\}$ and $p_{102} \in \left(-1,\frac12\right)$, this cubic has a degenerate global maximum at $e_1$ and two non-degenerate global maxima at the points \eqref{nondeg_zeros} on $S^2$.
\end{itemize} }

\begin{proof}
Expression \eqref{p_degenerate_maximum} is obtained from \eqref{e1face} by setting $p_{030} = p_{111} = 0$, $p_{120} = \frac12$. By construction the set ${\cal F}_4$ consists of exactly those cubics $p \in B_1(S^2)$ such that the restriction of $1-p$ on the great circle through $e_1,e_2$ has order of smallness at least 4 at $x = e_1$. By Lemma \ref{lem:order_face} it is a face of $B_1(S^2)$.

Functions \eqref{p12_p03_functions} simplify to
\begin{align*}
p_{03}(\varphi) &= \sin\varphi\left(p_{003}\sin^2\varphi + 3p_{012}\cos\varphi\sin\varphi + 3p_{021}\cos^2\varphi\right), \\
p_{12}(\varphi) &= \left(p_{102} - \frac12\right)\sin^2\varphi + \frac12.
\end{align*}
Hence by Lemma \ref{lem:face_S2} a polynomial $p$ of the form \eqref{p_degenerate_maximum} is in $B_1(S^2)$ if and only if the LMI
\[ \begin{pmatrix} 6\left(\frac12 - p_{102}\right) & 2\mu(\varphi) & -2\sin\varphi\left(\frac12 - p_{102}\right) \\ 2\mu(\varphi) & 3 - 2\sin^2\varphi\left(\frac12 - p_{102}\right) & 0 \\ -2\sin\varphi\left(\frac12 - p_{102}\right) & 0 & 1 \end{pmatrix} \succeq 0
\]
holds for every $\varphi \in [-\pi,\pi]$. Here we denoted $\mu(\varphi) = p_{003}\sin^2\varphi + 3p_{012}\cos\varphi\sin\varphi + 3p_{021}\cos^2\varphi$ for short.

Let us show that this condition defines a 4-dimensional convex compact body ${\cal B}_4$ in the space of coefficients $(p_{102},p_{021},p_{012},p_{003})$, and hence the face ${\cal F}_4$ of $B_1(S^2)$ is 4-dimensional.

By taking the Schur complement of the lower right corner element, we obtain the equivalent condition
\[ \begin{pmatrix} (1 - 2p_{102})\left(3 - \sin^2\varphi(1 - 2p_{102})\right) & 2\mu(\varphi) \\ 2\mu(\varphi) & 3 - \sin^2\varphi(1 - 2p_{102}) \end{pmatrix} \succeq 0 \quad \forall\ \varphi \in [-\pi,\pi],
\]
which simplifies to
\[ \pm2\mu(\varphi) \leq \sqrt{1 - 2p_{102}}\left(3\cos^2\varphi + 2\sin^2\varphi(1 + p_{102})\right) \quad \forall\ \varphi \in [-\pi,\pi].
\]
This can be written as 
\[ \pm 2 \begin{pmatrix} \cos\varphi \\ \sin\varphi \end{pmatrix}^T \begin{pmatrix} 3p_{021} & \frac32p_{012} \\ \frac32p_{012} & p_{003} \end{pmatrix} \begin{pmatrix} \cos\varphi \\ \sin\varphi \end{pmatrix} \preceq \sqrt{1 - 2p_{102}} \begin{pmatrix} \cos\varphi \\ \sin\varphi \end{pmatrix}^T \begin{pmatrix} 3 & 0 \\ 0 & 2(1 + p_{102}) \end{pmatrix} \begin{pmatrix} \cos\varphi \\ \sin\varphi \end{pmatrix}
\]
\[ \forall\ \varphi \in [-\pi,\pi]
\]
and simplifies to the matrix inequalities
\begin{equation} \label{degenerate_maxima_LMI}
\sqrt{1 - 2p_{102}}\begin{pmatrix} 3 & 0 \\ 0 & 2(1 + p_{102}) \end{pmatrix} \pm \begin{pmatrix} 6p_{021} & 3p_{012} \\ 3p_{012} & 2p_{003} \end{pmatrix} \succeq 0.
\end{equation}

Although these are not linear matrix inequalities, they are linear in the three variables $p_{021},p_{012},p_{003}$. Hence it is convenient to consider the 4-dimensional set ${\cal B}_4$ as a union of 3-dimensional convex compact slices, each of which is obtained by fixing the value $p_{102} \in \left[-1,\frac12\right]$ and analyzing LMIs \eqref{degenerate_maxima_LMI} in the remaining coefficients $(p_{021},p_{012},p_{003})$. A closer look reveals that each such slice is linearly equivalent to the intersection of two 3-dimensional Lorentz cones with tips pointing away from each other.

For a sample of values $p_{102} \in \left[-1,\frac12\right]$ the slices of ${\cal B}_4$ are visualized in Fig.~\ref{fig:sectionsB4}. The section of ${\cal B}_4$ is
\begin{itemize}
\item a line segment given by $p_{012} = p_{003} = 0$, $|p_{021}| \leq \frac{\sqrt{3}}{2}$ if $p_{102} = -1$;
\item the 3-dimensional centrally symmetric intersection of two ellipsoidal cones if $p_{102} \in \left(-1,\frac12\right)$;
\item the origin if $p_{102} = \frac12$,
\end{itemize}
and ${\cal B}_4$ is indeed 4-dimensional.

\begin{figure}
\centering
\includegraphics[width=13.62cm,height=7.18cm]{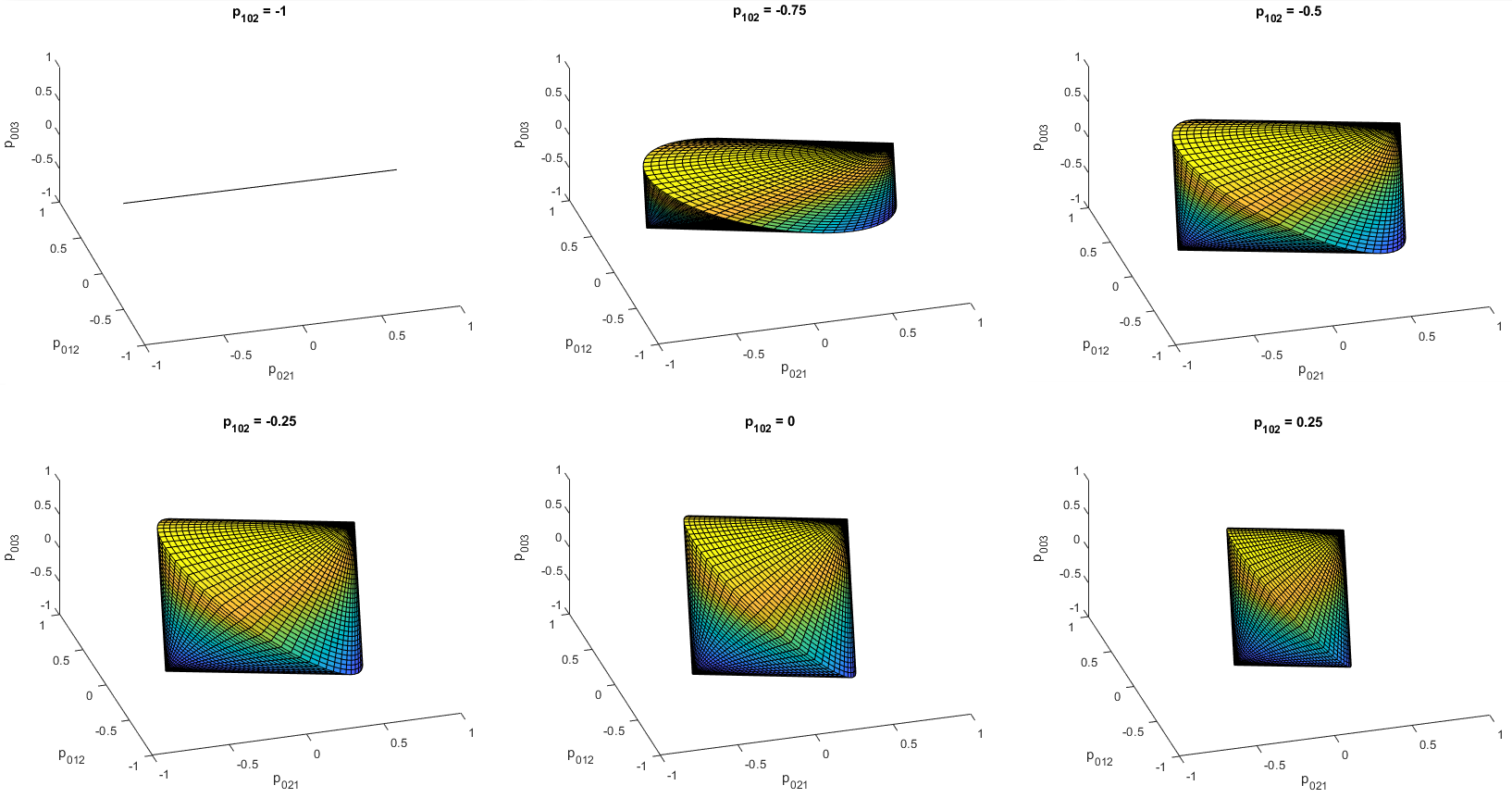}
\caption{Sections of the body ${\cal B}_4$ for different values of $p_{102}$. For $p_{102} \to \frac12$ the section shrinks to a point, for $p_{102} \to -1$ to a line segment.}
\label{fig:sectionsB4}
\end{figure}

We now look for the extremal points of ${\cal B}_4$, which must also be extremal in their corresponding section.

The end-points $\left(p_{102},p_{021},p_{012},p_{003}\right) = \left(-1,\pm\frac{\sqrt{3}}{2},0,0\right)$ of the line segment from the first item are extremal in ${\cal B}_4$ and correspond to the cubics $p(x) = x_1^3 + 3\left(\frac12x_2^2 - x_3^2\right)x_1 \pm \frac{3\sqrt{3}}{2}x_2^2x_3$. In this case we have $p_{12}(\pm\pi/2) = -1$, and these cubics have already been described in Corollary \ref{cor:3maxima_extremal}. Clearly they correspond to the first item in this corollary.

The point $\left(p_{102},p_{021},p_{012},p_{003}\right) = \left(\frac12,0,0,0\right)$ from the third item above is extremal in ${\cal B}_4$ and corresponds to the cubic $p(x) = x_1^3 + \frac32(x_2^2 + x_3^2)x_1$. Clearly this cubic is zonal and has a unique degenerate maximum on $S^2$ at $x = e_1$.

If $p_{102}$ assumes an intermediate value, then the extremal points of the corresponding section lie at the tips of the intersecting ellipsoidal cones and at the intersection of their boundaries.

Let us first consider the tips $\left(p_{021},p_{012},p_{003}\right) = \pm\sqrt{1-2p_{102}}\left(\frac12,0,(1+p_{102})\right)$. Since the two cubics corresponding to the tips are mapped to each other under the transformation $(x_1,x_2,x_3) \mapsto (x_1,x_2,-x_3)$ of $\mathbb R^3$, we shall consider only one of them, e.g.,
\[ p(x) = x_1^3 + \frac32\left(x_2^2 + 2p_{102}x_3^2\right)x_1 + \sqrt{1-2p_{102}}\left(\frac32x_2^2x_3 + (1+p_{102})x_3^3\right).
\]
Let us evaluate this cubic on the sphere $S^2$. Replacing $x_2^2$ by $1 - x_1^2 - x_3^2$, we obtain
\begin{align*}
p(x) &= x_1^3 + \frac32\left(1 - x_1^2 - (1-2p_{102})x_3^2\right)x_1 + \sqrt{1-2p_{102}}x_3\left(\frac32(1 - x_1^2) - \frac12(1-2p_{102})x_3^2\right) \\
&= -\frac12\left(x_1 + \sqrt{1-2p_{102}}x_3\right)^3 + \frac32\left(x_1 + \sqrt{1-2p_{102}}x_3\right).
\end{align*}
Hence the cubic $p$ is zonal, with the global maximum achieved on the intersection of $S^2$ with the plane given by the equation $x_1 + \sqrt{1-2p_{102}}x_3 = 1$. By Lemma \ref{lem:4_planar_maxima} it is an extremal point of $B_1(S^2)$.

We now pass to the intersection of the boundaries of the ellipsoidal cones defined by LMIs \eqref{degenerate_maxima_LMI}. In order to compute this intersection we have to equate the determinants of both the sum and the difference in \eqref{degenerate_maxima_LMI} to zero, yielding
\[ p_{003} + 2(1+p_{102})p_{021} = 0, \quad 3p_{012}^2 + 8(1+p_{102})p_{021}^2 = 2(1+p_{102})(1-2p_{102}).
\]
Let us parameterize the locus of solutions by an angle $\xi \in [-\pi,\pi]$, setting
\[ p_{003} = \sqrt{1-2p_{102}}(1+p_{102})\sin\xi,\quad p_{021} = -\frac12\sqrt{1-2p_{102}}\sin\xi, \quad p_{012} = \frac{\sqrt{6}}{3}\sqrt{1-2p_{102}}\sqrt{1+p_{102}}\cos\xi.
\]
Then the cubic $p(x)$ is given by
\begin{equation} \label{deg_nondeg_nondeg}
x_1^3 + 3\left(\frac12x_2^2 + p_{102}x_3^2\right)x_1 + \sqrt{1-2p_{102}}\left(-\frac32\sin\xi x_2^2x_3 + \sqrt{6}\sqrt{1+p_{102}}\cos\xi x_2x_3^2 + (1+p_{102})\sin\xi x_3^3\right).
\end{equation}
It is checked straightforwardly that this cubic assumes the value 1 at the two points
\begin{equation} \label{nondeg_zeros}
x^{\pm} = \frac{1}{2\left(2 - p_{102} \pm (1-2p_{102})\sin\xi\right)}\begin{pmatrix} 1 + 4p_{102} \mp (1-2p_{102})\sin\xi \\ \sqrt{6(1-2p_{102})(1+p_{102})}\cos\xi \\ 3\sqrt{1-2p_{102}}(\pm1+\sin\xi) \end{pmatrix} \in S^2.
\end{equation}
For fixed $p_{102}$ and variable $\xi$, the maxima $x^{\pm}$ run around the circles
\begin{equation} \label{circles_triply_degenerated}
C^{\pm} = S^2 \cap \left\{ x \left| x_1 \pm \sqrt{1-2p_{102}}x_3 = 1 \right. \right\},
\end{equation}
respectively.

If $\xi = \pm\frac{\pi}{2}$, then the point $x^{\mp}$ coincides with $e_1$. In this case the restriction of $p$ on the circle $C^{\mp}$ has order of smallness 6 at $x = e_1$. Indeed, setting $x_2^2 = 1 - x_1^2 - x_3^2$ and $x_1 = 1 \pm \sqrt{1-2p_{102}}x_3$ in the expression \eqref{deg_nondeg_nondeg} for $p$, we obtain $p(x) = 1 \pm 2\sqrt{1-2p_{102}}(1+p_{102})x_3^3$ on $C^{\mp}$. However, $x_3$ itself is already quadratic in the parameter of the circle at $x = e_1$. If the point $x^{\pm}$ does not coincide with $e_1$, then it is an isolated non-degenerate maximum of $p$ on $S^2$. Since the maxima determine all coefficients of $p$ uniquely, the cubic $p$ is an extremal point of $B_1(S^2)$.
\end{proof}

\subsection{Four distinct non-degenerate global maxima}

We now consider extremal cubics $p \in B_1(S^2)$ which have only non-degenerate global maxima. By Lemma \ref{lem:min4maxima} there must be at least 4 of them, and by Lemma \ref{lem:4_planar_maxima} they have to be affinely independent, i.e., not lying in the same plane. An explicit parametrization of the coefficients of such extremal cubics turns out to be too complicated, and therefore we shall describe these cubics by means of the location of their maxima on $S^2$.

It turns out that critical points of a homogeneous cubic cannot be scattered arbitrarily over the sphere. However, if a quadruple of points $u^1,u^2,u^3,u^4 \in S^2$ is critical for some cubic, then the image of the quadruple under an orthogonal transformation is also critical for some cubic. Hence it is sufficient to investigate the Gramian $\Gamma$ of the critical points. The following result provides a necessary condition on this Gramian.

{\lemma \label{lem:4maxima_necessary} Let $p$ be a homogeneous cubic on $\mathbb R^3$. Suppose $p$ has four affinely independent critical points $u^1,u^2,u^3,u^4$ on $S^2$, such that no three of them lie on a great circle, and the values $p(u^i)$ equal 1 for all $i = 1,2,3,4$. Then there exists a vector $z \in \mathbb R^4$ with non-zero entries such that $\sum_{i=1}^4 z_i = 1$ and the Gramian of the critical points $u^i$ is given by
\begin{equation} \label{Gamma_z}
\Gamma = \begin{pmatrix} 1 & -\frac12+\frac16g_1g_2 & -\frac12+\frac16g_1g_3 & -\frac12+\frac16g_1g_4 \\ -\frac12+\frac16g_1g_2 & 1 & -\frac12+\frac16g_2g_3 & -\frac12+\frac16g_2g_4 \\ -\frac12+\frac16g_1g_3 & -\frac12+\frac16g_2g_3 & 1 & -\frac12+\frac16g_3g_4 \\ -\frac12+\frac16g_1g_4 & -\frac12+\frac16g_2g_4 & -\frac12+\frac16g_3g_4 & 1 \end{pmatrix},
\end{equation}
where $g_i = -3+z_i^{-1}$, $i = 1,2,3,4$.

On the other hand, let $u^1,u^2,u^3,u^4 \in S^2$ be affinely independent points such that their Gramian has the form \eqref{Gamma_z} for some $z \in \mathbb R^4$ with non-zero components and satisfying $\sum_{i=1}^4 z_i = 1$. Then there exists a unique homogeneous cubic $p$ such that the $u^i$ are critical points of $p$ on $S^2$ and $p(u^i) = 1$ for all $i = 1,2,3,4$. }

\begin{proof}
Let $p$ be as in the first part of the lemma, and let $\Gamma$ be the Gramian of the critical points $u^i$. Since the $u^i$ are affinely independent, they can be used as an affine basis defining barycentric coordinates in $\mathbb R^3$. In other words, for every $x \in \mathbb R^3$ there exists a unique vector $b \in \mathbb R^4$ such that $\sum_{i=1}^4 b_i = 1$ and $x = \sum_{i=1}^4 b_iu^i$. Let $z$ be the barycentric coordinate vector of the origin in $\mathbb R^3$.

If at least one of the components of $z$ vanishes, then the origin is an affine combination of at most three of the critical points $u^i$, and these critical points must lie on a great circle, in contradiction with the assumptions of the lemma. Therefore the components of $z$ are non-zero, and the quantities $g_i$ are well-defined.

Since all critical points $u^i$ are located on the unit sphere, the diagonal elements of the Gramian equal 1. 
Let us show that its off-diagonal elements $\Gamma_{ij}$ have the values claimed in \eqref{Gamma_z}.

Define a linear projection $\Pi: \mathbb R^4 \to \mathbb R^3$ by the formula $y \mapsto \sum_{i=1}^4 y_iu^i$. Then $\Pi(z) = 0$ by definition of $z$, and the unit basis vectors ${\bf e}_i$ in $\mathbb R^4$ are mapped by $\Pi$ to the critical points $u^i$, $i = 1,2,3,4$.

Define a homogeneous cubic polynomial $q$ and a quadratic form $\Omega$ on $\mathbb R^4$ by the composition $q = p \circ \Pi$ and the formula $\Omega(y) = \|\Pi(y)\|_2^2 = \langle \Pi(y),\Pi(y) \rangle$, respectively. Let $q_{ijk},\Omega_{ij}$ be the coefficients of $q,\Omega$, respectively, i.e., $q_{ijk},\Omega_{ij}$ are symmetric in their indices and $q(y) = \sum_{i,j,k = 1}^4 q_{ijk}y_iy_jy_k$, $\Omega(y) = \sum_{i,j=1}^4 \Omega_{ij}y_iy_j$ for all $y \in \mathbb R^4$. Then we have $\Omega_{ij} = \langle \Pi({\bf e}_i),\Pi({\bf e}_j) \rangle = \langle u^i,u^j \rangle$, and the coefficient matrix of $\Omega$ is identical with the Gramian $\Gamma$ of the critical points $u^i$.

The gradient of $p$ at the critical point $u^i$ is orthogonal to $S^2$ and hence proportional to $u^i$. By virtue of homogeneity of $p$ we have $\nabla p(u^i) = 3u^i$. It follows that the gradient of the 0-homogeneous ratio $\frac{p(x)}{\|x\|_2^3}$ vanishes at $u^i$. But then the 0-homogeneous ratio $r(y) = \frac{q(y)}{\Omega(y)^{3/2}} = \left. \frac{p(x)}{\|x\|_2^3}\right|_{x = \Pi(y)}$, defined everywhere on $\mathbb R^4$ except $\Pi^{-1}[0]$, has vanishing gradient at ${\bf e}_i$. In particular, $r({\bf e}_i) = \frac{q_{iii}}{\Omega_{ii}^{3/2}} = q_{iii} = 1$, and by the rule of derivation of a quotient we get for every $i,j = 1,2,3,4$ that
\[ \left. \frac{\partial r(y)}{\partial y_j} \right|_{y = {\bf e}_i} = \frac{(3q_{iij}) \cdot \Omega_{ii}^{3/2} - q_{iii} \cdot \left( \frac32\Omega_{ii}^{1/2}\cdot 2\Omega_{ij} \right)}{\Omega_{ii}^3} = 3(q_{iij} - \Gamma_{ij}) = 0.
\]
It follows that $q_{iij} = q_{ijj} = \Gamma_{ij}$, $i,j = 1,2,3,4$.

Finally, for every $\alpha \in \mathbb R$ and every $y \in \mathbb R^4$ we have $\Pi(y + \alpha z) = \Pi(y)$ and hence
\[ \langle \nabla q(y),z \rangle = 3\sum_{i,j,k=1}^4q_{ijk}y_iy_jz_k = 0\quad \forall\ y \in \mathbb R^4.
\]
This yields 10 linearly independent relations $\sum_{k=1}^4q_{ijk}z_k = 0$, $i,j = 1,\dots,4$ on the coefficients $q_{ijk}$. Writing this system out, we obtain
\[ \begin{pmatrix}
z_1+z_2 & & & & & & & & z_4 & z_3 \\
& z_1+z_3 & & & & & & z_4 & & z_2 \\
& & z_1+z_4 & & & & & z_3 & z_2 & \\
& & & z_2+z_3 & & & z_4 & & & z_1 \\
& & & & z_2+z_4 & & z_3 & & z_1 & \\
& & & & & z_3+z_4 & z_2 & z_1 & & \\
z_2 & z_3 & z_4 & & & & & & & \\
z_1 & & & z_3 & z_4 & & & & & \\
& z_1 & & z_2 & & z_4 & & & & \\
& & z_1 & & z_2 & z_3 & & & &
\end{pmatrix} \begin{pmatrix} \Gamma_{12} \\ \Gamma_{13} \\ \Gamma_{14} \\ \Gamma_{23} \\ \Gamma_{24} \\ \Gamma_{34} \\ q_{234} \\ q_{134} \\ q_{124} \\ q_{123} \end{pmatrix} = -\begin{pmatrix} 0 \\ 0 \\ 0 \\ 0 \\ 0 \\ 0 \\ z_1 \\ z_2 \\ z_3 \\ z_4 \end{pmatrix}.
\]
The determinant of the resulting $10 \times 10$ matrix equals $-12z_1^2z_2^2z_3^2z_4^2(z_1 + z_2 + z_3 + z_4)^2 \not= 0$, and the off-diagonal elements of the Gramian as well as the totally mixed coefficients $q_{ijk}$ can be computed as explicit functions of $z$ by inverting this matrix. This yields \eqref{Gamma_z} and the expressions
\begin{align*}
q_{123} &= 1 + \frac{1 - 2z_4 - 3z_1z_2 - 3z_1z_3 - 3z_2z_3}{6z_1z_2z_3} = 1 - \frac{g_1g_2 + g_1g_3 + g_2g_3 + g_1g_2g_3}{6}, \\
q_{124} &= 1 + \frac{1 - 2z_3 - 3z_1z_2 - 3z_1z_4 - 3z_2z_4}{6z_1z_2z_4} = 1 - \frac{g_1g_2 + g_1g_4 + g_2g_4 + g_1g_2g_4}{6}, \\
q_{134} &= 1 + \frac{1 - 2z_2 - 3z_1z_3 - 3z_1z_4 - 3z_3z_4}{6z_1z_3z_4} = 1 - \frac{g_1g_3 + g_1g_4 + g_3g_4 + g_1g_3g_4}{6}, \\
q_{234} &= 1 + \frac{1 - 2z_1 - 3z_2z_3 - 3z_2z_4 - 3z_3z_4}{6z_2z_3z_4} = 1 - \frac{g_2g_3 + g_2g_4 + g_3g_4 + g_2g_3g_4}{6}
\end{align*}
for the remaining coefficients.

\medskip

Let us now prove the second part of the lemma. It is verified by direct calculus that $z$ is a kernel vector of the matrix $\Gamma$ in \eqref{Gamma_z}. Hence $\sum_{i=1}^4 z_iu^i = 0$, and $z$ is the barycentric coordinate of the origin $0 \in \mathbb R^3$ with respect to the affine basis $\{u^i\}_{i = 1,2,3,4}$. Let $\Pi$ be the linear projection defined above. Define a homogeneous cubic $q$ and a quadratic form $\Omega$ on $\mathbb R^4$ by the coefficient functions $q_{ijk},\Omega_{ij}$ depending explicitly on $z$ as above. Then the basis vectors ${\bf e}_i$ are critical points of the ratio $r(y) = \frac{q(y)}{\Omega(y)^{3/2}}$, $\Omega(y) \equiv \|\Pi(y)\|_2^2$, and $\langle \nabla q(y),z \rangle \equiv 0$ for all $y \in \mathbb R^4$. Define the homogeneous cubic $p$ by the relation $p(x) = q(b(x))$, where $b(x)$ is the barycentric coordinate of $x$. Then $q = p \circ \Pi$, and by reversing the reasoning above we obtain that $u^i = \Pi({\bf e}_i)$ are critical points of $p$ on $S^2$. Finally, performing the chain of arguments in the first part of the proof with an arbitrary homogeneous cubic $\tilde p$ having critical points at $u^i$ with values $p(u^i) = 1$ leads to the same polynomial $q$, and hence $\tilde p$ coincides with $p$.
\end{proof}

Since the location of the quadruple of critical points determines the cubic $p$ uniquely, it also determines the type of the critical points. Recall that if the critical points are global maxima, then Lemma \ref{lem:max_angle} yields an upper bound $\frac{2\pi}{3}$ on the angles between them. Equivalently, the off-diagonal elements of the Gramian are lower-bounded by $-\frac12$. We now show that this simple condition is essentially sufficient to guarantee that the critical points are indeed maxima.

{\lemma \label{lem:angle_maxima} Assume the conditions of Lemma \ref{lem:4maxima_necessary}. If the off-diagonal elements of the Gramian $\Gamma$ lie in the interval $\left(-\frac12,1\right)$, then the critical points $u^i$ are all non-degenerate local maxima of $p$ on $S^2$. In this case either $z_i \in \left(0,\frac13\right)$ for all $i$, or $z_i > \frac13$ for three indices $i$, and the remaining component of $z$ is negative. }

\begin{proof}
Consider again the 0-homogeneous function $r(y)$ introduced in the proof of the previous lemma. Its Hessian at the basis vector ${\bf e}_i$ has rank at most 2, because $r \equiv 1$ on the linear subspace spanned by ${\bf e}_i$ and $z$. The corresponding two non-trivial eigenvalues of this Hessian then determine the type of the critical point $u^i$.

Without loss of generality, consider the basis vector ${\bf e}_1$. The Hessian of $r$ at this point evaluates to
\[ \begin{pmatrix} 0 & 0 & 0 & 0 \\
0 &  -\frac{(g_1g_2 - 9)^2}{12} & \frac{3z_2 - 6z_1 + 3z_3 + 9(z_1z_2 + z_1z_3 - z_2z_3 + z_1^2) - 1}{12z_1^2z_2z_3} & \frac{3z_2 - 6z_1 + 3z_4 + 9(z_1z_2 + z_1z_4 - z_2z_4 + z_1^2) - 1}{12z_1^2z_2z_4} \\
0 & \frac{3z_2 - 6z_1 + 3z_3 + 9(z_1z_2 + z_1z_3 - z_2z_3 + z_1^2) - 1}{12z_1^2z_2z_3} & -\frac{(g_1g_3 - 9)^2}{12} & \frac{3z_3 - 6z_1 + 3z_4 + 9(z_1z_3 + z_1z_4 - z_3z_4 + z_1^2) - 1}{12z_1^2z_3z_4} \\
0 & \frac{3z_2 - 6z_1 + 3z_4 + 9(z_1z_2 + z_1z_4 - z_2z_4 + z_1^2) - 1}{12z_1^2z_2z_4} & \frac{3z_3 - 6z_1 + 3z_4 + 9(z_1z_3 + z_1z_4 - z_3z_4 + z_1^2) - 1}{12z_1^2z_3z_4} & -\frac{(g_1g_4 - 9)^2}{12} \end{pmatrix}.
\]
Since $\Gamma_{ij} < 1$ for $i \not= j$, we have $g_ig_j < 9$, and all diagonal elements of the Hessian except the first one are strictly negative. Its non-trivial principal minors of order two evaluate to
\[ \frac{g_2g_3g_4z_2}{12z_1^3z_3z_4}, \quad \frac{g_2g_3g_4z_3}{12z_1^3z_2z_4}, \quad \frac{g_2g_3g_4z_4}{12z_1^3z_2z_3},
\]
respectively. The condition $\Gamma_{ij} > -\frac12$ implies that all $g_i$ are non-zero and have the same sign. Hence also all non-trivial second order principal minors are non-zero have the same sign, namely that of the product $g_1z_1z_2z_3z_4$. Let us show that this product is positive.

If all $g_i$ are all positive, then $z_i \in \left(0,\frac13\right)$ for all $i$, and the claim follows.

If all $g_i$ are negative, then $z_i \in (-\infty,0) \cup \left(\frac13,+\infty\right)$ for all $i$. Since the $z_i$ sum to 1, they cannot be all positive, otherwise $\sum_{i=1}^4 z_i > \frac43$. On the other hand, if two distinct components $z_i,z_j$ are negative, then $g_i,g_j < -3$ and $g_ig_j > 9$, in contradiction with the condition $\Gamma_{ij} < 1$. Thus exactly one of the $z_i$ is negative, and the other components of $z$ are strictly larger than $\frac13$. In addition, $g_1z_1z_2z_3z_4 > 0$ also in this case.

Hence the Hessian is negative semi-definite of rank 2, and $u^1$ is a non-degenerate (local) maximum of $p$ on $S^2$. For the other critical points the claim is proven similarly.
\end{proof}

In the first case presented in Lemma \ref{lem:angle_maxima}, the conditions $z_i \in \left(0,\frac13\right)$ and $\sum_{i=1}^4 z_i = 1$ define an open 3-dimensional simplex in the space of variables $z \in \mathbb R^4$, which can be parameterized by a barycentric coordinate $b \in \Delta_4^o$. In the second case the conditions $z_1,z_2,z_3 > \frac13$, $\sum_{i=1}^4 z_i = 1$ allow to parameterize $z$ by a variable $b \in \mathbb R_{++}^3$. This yields the following result.

{\corollary \label{cor:4global_max_necessary} Let $p \in \partial B_1(S^2)$ have four distinct non-degenerate affinely independent global maxima $u^1,u^2,u^3,u^4$ on $S^2$, such that no three of them lie on a great circle. Then exactly one of the following conditions holds:
\begin{itemize}
\item there exists $b \in \Delta_4^o$ such that the Gramian $\Gamma$ of the maxima $u^i$ has the form \eqref{Gamma_z} with $z_i = \frac{1-b_i}{3}$ or equivalently, $g_i = \frac{3}{b_i^{-1}-1}$, $i = 1,2,3,4$;
\item there exists a permutation matrix $P$ and numbers $b_1,b_2,b_3 > 0$ such that the product $P\Gamma P^T$ has the form \eqref{Gamma_z} with $z_i = \frac{1+b_i}{3}$, $i = 1,2,3$, $z_4 = -\frac{b_1+b_2+b_3}{3}$ or equivalently, $g_i = -\frac{3}{b_i^{-1}+1}$, $i = 1,2,3$, $g_4 = -\frac{3(b_1+b_2+b_3+1)}{b_1+b_2+b_3}$.
\end{itemize}
In both cases $p$ is an extremal point of $B_1(S^2)$. }

\begin{proof}
Since $p \in \partial B_1(S^2)$, its maximum over $S^2$ equals 1. Hence $p$ satisfies the conditions of the first part of Lemma \ref{lem:4maxima_necessary}, and $\Gamma$ is of the form \eqref{Gamma_z} for some vector $z$ with non-zero entries and $\sum_{i=1}^4 z_i = 1$. By Lemma \ref{lem:max_angle} the angles between the maxima $u^i$ are strictly smaller than $\frac{2\pi}{3}$, and hence the off-diagonal elements of the Gramian $\Gamma$ lie in the interval $\left(-\frac12,1\right)$. By Lemma \ref{lem:angle_maxima} the elements of $z$ can be represented in the form stated in the corollary. Finally, by the second part of Lemma \ref{lem:4maxima_necessary} $p$ is determined uniquely by the quadruple $u^i$, and hence it is an extremal point of $B_1(S^2)$.
\end{proof}

We still do not know whether every choice of $b$ in the two cases of Corollary \ref{cor:4global_max_necessary} leads to an orbit of extremal points of $B_1(S^2)$. The corresponding matrix \eqref{Gamma_z} may not be a Gramian of four points $u^i$ on $S^2$, and even if it is, the polynomial $p$ defined by such a quadruple of points is guaranteed to have only \emph{local} maxima at $u^i$ and thus may lie outside of $B_1(S^2)$. First we shall show that in either case every choice of $b$ leads to a valid Gramian.

{\lemma \label{lem:4maxima_centre} Let $b = (b_1,b_2,b_3,b_4)^T \in \Delta_4^o$ be a point in the interior of the simplex in $\mathbb R^4$, and define the matrix
\begin{equation} \label{Gamma_central}
\Gamma = \begin{pmatrix} 1 & -\frac12+\frac32\frac{1}{(b_1^{-1}-1)(b_2^{-1}-1)} & -\frac12+\frac32\frac{1}{(b_1^{-1}-1)(b_3^{-1}-1)} & -\frac12+\frac32\frac{1}{(b_1^{-1}-1)(b_4^{-1}-1)} \\
-\frac12+\frac32\frac{1}{(b_1^{-1}-1)(b_2^{-1}-1)} & 1 & -\frac12+\frac32\frac{1}{(b_2^{-1}-1)(b_3^{-1}-1)} & -\frac12+\frac32\frac{1}{(b_2^{-1}-1)(b_4^{-1}-1)} \\
-\frac12+\frac32\frac{1}{(b_1^{-1}-1)(b_3^{-1}-1)} & -\frac12+\frac32\frac{1}{(b_2^{-1}-1)(b_3^{-1}-1)} & 1 & -\frac12+\frac32\frac{1}{(b_3^{-1}-1)(b_4^{-1}-1)} \\
-\frac12+\frac32\frac{1}{(b_1^{-1}-1)(b_4^{-1}-1)} & -\frac12+\frac32\frac{1}{(b_2^{-1}-1)(b_4^{-1}-1)} & -\frac12+\frac32\frac{1}{(b_3^{-1}-1)(b_4^{-1}-1)} & 1 \end{pmatrix}.
\end{equation}
Then $\Gamma$ is of the form \eqref{Gamma_z} with $z_i = \frac{1-b_i}{3}$ or equivalently, $g_i = \frac{3}{b_i^{-1}-1}$, $i = 1,2,3,4$, it is positive semi-definite of rank 3, and its off-diagonal elements are contained in the interval $(-\frac12,1)$. The vector $z = (z_1,\dots,z_4)^T$ is in $\Delta_4^o$. }

\begin{proof}
We have $b_i^{-1}-1 > 0$ for all $i$, and hence $\Gamma_{ij} > -\frac12$. Further, for $i \not= j$ we have $b_i + b_j < 1$, implying $(b_i^{-1}-1)(b_j^{-1}-1) > 1$ and $\Gamma_{ij} < 1$. It also follows that all principal submatrices of $\Gamma$ of size 1 or 2 have a positive determinant.

The determinant of the $3 \times 3$ principal upper left submatrix of $\Gamma$ equals $\frac{9b_1b_2b_3b_4(b_1^{-1}+b_2^{-1}+b_3^{-1}+b_4^{-1}-4)}{4(1-b_1)^2(1-b_2)^2(1-b_3)^2} > 0$, similarly the other principal minors of order 3 are positive. It is verified directly that $\det\,\Gamma = 0$. Hence $\Gamma$ is positive semi-definite of rank 3 by Sylvester's criterion.

Finally, it is checked straightforwardly that $\Gamma$ has the form \eqref{Gamma_z} with $z \in \Delta_4^o$ given by the claimed expression.
\end{proof}

{\lemma \label{lem:4maxima_wing} Let $b_1,b_2,b_3$ be positive real numbers, and define the matrix
\begin{equation} \label{Gamma_wing}
\Gamma = \begin{pmatrix} 1 & -\frac12+\frac32\frac{1}{(b_1^{-1}+1)(b_2^{-1}+1)} & -\frac12+\frac32\frac{1}{(b_1^{-1}+1)(b_3^{-1}+1)} & -\frac12+\frac32\frac{(b_1+b_2+b_3+1)b_1}{(b_1+b_2+b_3)(b_1+1)} \\
-\frac12+\frac32\frac{1}{(b_1^{-1}+1)(b_2^{-1}+1)} & 1 & -\frac12+\frac32\frac{1}{(b_2^{-1}+1)(b_3^{-1}+1)} & -\frac12+\frac32\frac{(b_1+b_2+b_3+1)b_2}{(b_1+b_2+b_3)(b_2+1)} \\
-\frac12+\frac32\frac{1}{(b_1^{-1}+1)(b_3^{-1}+1)} & -\frac12+\frac32\frac{1}{(b_2^{-1}+1)(b_3^{-1}+1)} & 1 & -\frac12+\frac32\frac{(b_1+b_2+b_3+1)b_3}{(b_1+b_2+b_3)(b_3+1)} \\
-\frac12+\frac32\frac{(b_1+b_2+b_3+1)b_1}{(b_1+b_2+b_3)(b_1+1)} & -\frac12+\frac32\frac{(b_1+b_2+b_3+1)b_2}{(b_1+b_2+b_3)(b_2+1)} & -\frac12+\frac32\frac{(b_1+b_2+b_3+1)b_3}{(b_1+b_2+b_3)(b_3+1)} & 1 \end{pmatrix}.
\end{equation}
Then $\Gamma$ is of the form \eqref{Gamma_z} with $z_i = \frac{1+b_i}{3}$, $i = 1,2,3$, $z_4 = -\frac{b_1+b_2+b_3}{3}$ or equivalently, $g_i = -\frac{3}{b_i^{-1}+1}$, $i = 1,2,3$, $g_4 = -\frac{3(b_1+b_2+b_3+1)}{b_1+b_2+b_3}$. The matrix $\Gamma$ is positive semi-definite of rank 3, and its off-diagonal elements are contained in the interval $(-\frac12,1)$. The vector $z = (z_1,\dots,z_4)^T$ has three positive and one negative entry, and $z_1+\dots+z_4 = 1$. }

\begin{proof}
Note that if $0 < \alpha < \beta$, then $0 < \frac{(\beta+1)\alpha}{\beta(\alpha+1)} < 1$. It follows that all off-diagonal elements of $\Gamma$ are contained in the interval $(-\frac12,1)$, and all its principal minors of order 1 or 2 are positive.

The upper left principal $3 \times 3$ submatrix of $\Gamma$ has determinant
\[ \frac{9(4b_1b_2b_3(b_1 + b_2 + b_3 + \frac34) + (b_1 + b_2 + b_3 + 1)(b_1b_2 + b_1b_3 + b_2b_3))}{4(b_1 + 1)^2(b_2 + 1)^2(b_3 + 1)^2} > 0.
\]
Therefore this submatrix is positive definite by Sylvester's criterion. By direct calculus it is verified that the matrix $\Gamma$ itself is singular, and consequently must be positive semi-definite of rank 3.

It is checked straightforwardly that $\Gamma$ has the form \eqref{Gamma_z} with $z$ having the claimed properties.
\end{proof}

The analysis is concluded by the following result.

{\lemma \label{lem:local_to_global} Let $p$ be a homogeneous cubic which has non-degenerate local maxima at the affinely independent points $u^1,u^2,u^3,u^4$ on $S^2$, such that the values $p(u^i)$ equal 1 for all $i = 1,2,3,4$. Then these maxima are global, and $p \in \partial B_1(S^2)$. }

\begin{proof}
To prove that $u^1,u^2,u^3,u^4$ are global maxima, it is sufficient to show that there are no other maxima on the sphere. We start by finding all other stationary points.

Consider the homogeneous cubic polynomial $q(y)$ and the quadratic form $\Omega(y)$  introduced in the proof of Lemma \ref{lem:4maxima_necessary}. To find the other stationary points of the polynomial, we write the first-order condition: there should exist $k \, \in \, \mathbb{R}$ such that $\nabla  q(y) = k \nabla \Omega(y)$, i.e. we have the following system of equations
\begin{eqnarray}
\label{eq:sys_opt}
    \frac{\partial q(y) }{\partial y_1} = k \frac{\partial \Omega(y)}{\partial y_1}, \quad
    \frac{\partial q(y) }{\partial y_2} = k \frac{\partial \Omega(y)}{\partial y_2}, \quad
    \frac{\partial q(y) }{\partial y_3} = k \frac{\partial \Omega(y)}{\partial y_3}, \quad
    \frac{\partial q(y) }{\partial y_4} = k \frac{\partial \Omega(y)}{\partial y_4}.
\end{eqnarray}
Note that $\langle z, \nabla  q(y) \rangle =  \langle z, \nabla \Omega(y) \rangle = 0$, so the equations are linearly dependent, and we need to consider only the first three equations. Since the system is invariant with respect to the replacement of $y \, \mapsto \, y + \alpha z$, we can assume $\sum_{i=1}^4 y_i = 0$ without loss of generality. 
Substituting $y_4 = -y_1 - y_2 - y_3$ into the first three equations in \eqref{eq:sys_opt}, we obtain a system of three equations in three variables $y_1, y_2, y_3$. Eliminating the scalar $k$, we get the following system of two equations 
\begin{equation}
\label{eq:eq_reduce_c}
    \frac{\partial q(y) }{\partial y_1}\cdot \frac{\partial \Omega(y)}{\partial y_2} - 
    \frac{\partial  q(y) }{\partial y_2} \cdot  \frac{\partial \Omega(y)}{\partial y_1} = 0
    , \quad \frac{\partial q(y) }{\partial y_1}\cdot \frac{\partial \Omega(y)}{\partial y_3} - 
    \frac{\partial  q(y)}{\partial y_3} \cdot  \frac{\partial \Omega(y)}{\partial y_1} = 0.
\end{equation}

Note that the l.h.s.  of the equations are homogeneous cubic polynomials in $y_1, y_2, y_3$. We consider these polynomials as cubic polynomials in $y_3$, i.e.
\begin{equation}
\label{eq:x3_equations}
    a_1 y_3^3 + b_1 y_3^2 + c_1 y_3 + d_1 = 0, \quad a_2 y_3^3 + b_2 y_3^2 + c_2 y_3 + d_2 = 0,
\end{equation}
where $a_i, b_i, c_i$ depend on $y_1, y_2$. Next, we find the resultant of the system of equations \eqref{eq:x3_equations}, i.e. the determinant of the matrix 
\begin{equation*}
 \text{Res} =   \text{det}\begin{pmatrix}
a_1   &  b_1 & c_1 &  d_1 & 0 &  0 \\
 0 &  a_1   &  b_1 & c_1 &  d_1  &  0 \\
 0 &  0 &  a_1   &  b_1 & c_1 &  d_1  \\
a_2   &  b_2 & c_2 &  d_2 & 0 &  0 \\
 0 &  a_2   &  b_2 & c_2 &  d_2  &  0 \\
 0 &  0 &  a_2   &  b_2 & c_2 &  d_2 
\end{pmatrix}.
\end{equation*}
Then, we calculate an analytic expression for the resultant and get the factorization
\begin{equation}
\label{resultant}
    \text{Res} = J(y_1, y_2)\cdot L(y_1, y_2)\cdot Q(y_1, y_2),
\end{equation}
where $L(y_1, y_2)$ is a homogeneous quadratic polynomial, $Q(y_1, y_2)$ a homogeneous cubic polynomial, and $J((y_1, y_2)$ a homogeneous polynomial of degree four.

Since we need to find a common solution for both equations \eqref{eq:x3_equations}, their resultant should equal zero. From the factorization \eqref{resultant} we obtain three sets of solutions. Note that we get each solution as a straight line in the hyperplane given by $\sum_{i=1}^4 y_i = 0$, corresponding to two opposite critical points of $p$ on the unit sphere $S^2$. 

The solutions of $J(y_1, y_2)=0$ correspond to borderline cases $z_i = \tfrac{1}{3}$ and projections of the maxima of $q$ onto hyperplane $\sum_{i=1}^4 y_i = 0$, i.e. $(1 - z_1, -z_2, -z_3,- z_4), \, \, ( - z_1, 1 -z_2, -z_3,- z_4), \, \, (- z_1, -z_2, 1 -z_3,- z_4), \, \,  ( - z_1, -z_2, -z_3,1 - z_4)$. 

If $L(y_1, y_2) = 0$ has solutions, then for them $\frac{\partial  q(x) }{\partial y_1} = 0 $ and $ \frac{\partial \Omega(x)}{\partial y_1} = 0$, and albeit being solutions of \eqref{eq:eq_reduce_c}, they do not correspond to solutions of \eqref{eq:sys_opt}. 

Finally, in the case $Q(y_1, y_2) = 0$  we get $3$ solutions of \eqref{eq:x3_equations}. Each solution corresponds to two opposite critical points on the sphere $S^2$. So, we get $6$ additional critical points on the sphere.

Our next step is to show that these $6$ solutions must be saddle points and can not be maxima. We employ a topological argument. We will show that if  $f\, : \, S^2 \rightarrow \R$ is an analytical function with  $4$ minima $x_i, 
\, i = 1,2,3, 4$ and $4$ maxima $y_j, 
\, j =1,2,3, 4$ on a sphere, such that
\begin{equation*}
    f(x_1) \leq  \ldots \leq f(x_4) \leq f(y_1)  \leq \ldots \leq f(y_4),
\end{equation*} then this function should have at least $6$ saddle points on the sphere.

Since $f$ is an analytic function, it reaches its minimum and maximum on the sphere. Then the range of the function is the segment and $f\, : \, S^2 \rightarrow I \subset \R$, there $I = [a, b].$ The function has a finite number of extrema, and they are all contained in $I$.

Consider a value $\gamma$ such that there is no critical point $\hat x$ such that $f(\hat x) = \gamma.$ Then we define the level set
\begin{equation*}
    X_{\gamma} = \left\{ x \, \in \, S^2\, | \, f(x) = \gamma \right\}.
\end{equation*}
This set does not contain critical points, and therefore the gradient of the function $f$ on this set is non-zero. Then by the implicit function theorem, the contour line is an analytic curve that extends through every point. Note that it cannot intersect itself because otherwise, the gradient is zero at the intersection point.  Moreover, it cannot be of infinite length, because otherwise there is an accumulation point, which must be critical. As a result, we obtain that $X_{\gamma}$ must be a closed curve with finite length.  Since the sphere $S^2$ is simply connected, any embedding of $S^1$ in $S^2$ is homotopic to any other embedding.  Thereby, the level line consists of a finite number of non-intersecting loops.

Let us for given $\gamma$ color the sphere such that all regions there $f(z) < \gamma$ are red and otherwise they are blue. In particular, for $\gamma < f(x_1)$ the whole sphere is red. Then we start to increase $\gamma$ and study the evolution of the coloring of $S^2$. After passing through a minimal point $x_i, \, i=1,2,3,4$ a new blue region around this point appears. 

Note that if we pass through a saddle point, there are two possibilities. 
\begin{figure}[!ht]
  \centering
  \begin{subfigure}[b]{0.4\linewidth}
    \includegraphics[width=\linewidth]{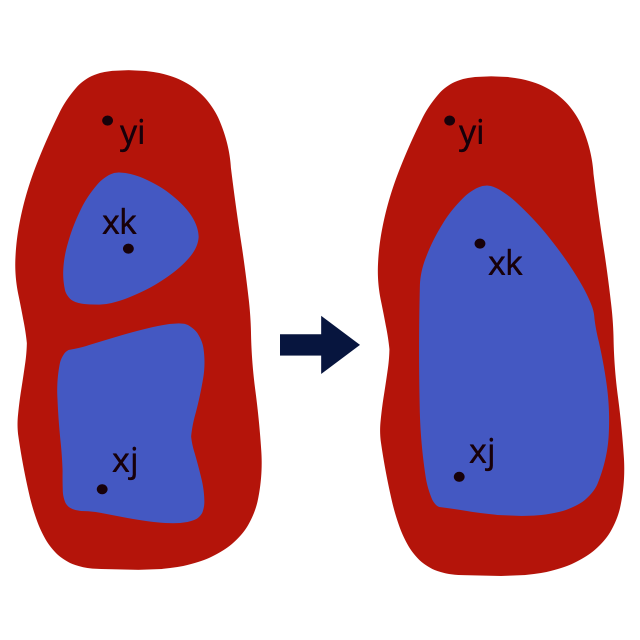}
     \caption{Two loops become one}
     \label{gr:3}
  \end{subfigure}
  \begin{subfigure}[b]{0.4\linewidth}
    \includegraphics[width=\linewidth]{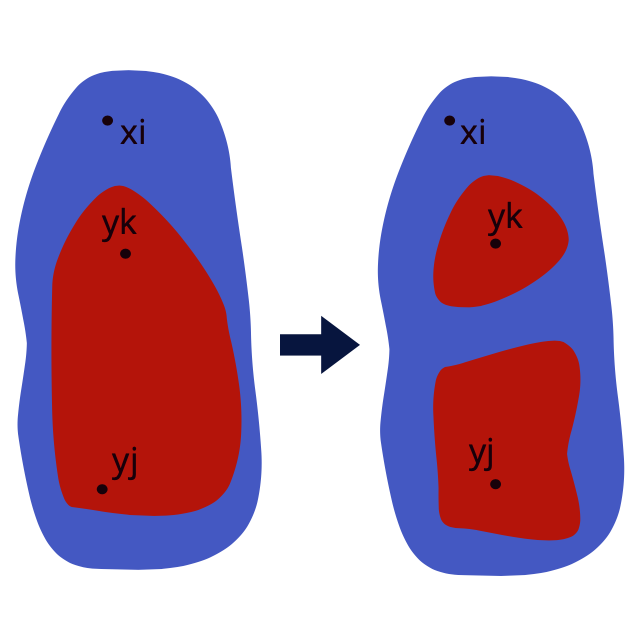}
    \caption{One loop splits into two}
    \label{gr:4}
  \end{subfigure}
  \caption{Regions representation after passing through a saddle-point}
\end{figure}
In the first case, two blue regions merge into one (two loops become one (Fig.\ref{gr:3})). Alternatively, one red region splits into two (one loop splits into two loops (Fig.\ref{gr:4})). 

Note that for $\gamma < f(y_i)$, we have a small red region around the maximum $y_i$, which will disappear after passing through this maximum.  And finally, for $\gamma > f(y_4)$, the whole sphere becomes blue. 
Then, we get that to connect all four blue regions, we need to pass through three saddle-points, and to split one red region onto four regions, we also need three saddle-points.

We obtain that for the case of $4$ minima and $4$ maxima, there must be at least $6$ saddle points.

So, since we already have four maximal and four minimal points and only six additional critical points, they can be only saddle points, and we do not have other maximal points. Moreover, since in our case the functional values at the maximal points $u^{i}, \, i=1,2,3,4$ are the same and equal 1, they are precisely the global maxima.
\end{proof}

\subsection{Final result}

We obtain the following classification of extremal points of $B_1(S^2)$, with the types of isolated global maxima described in Definition \ref{def:maxima_types}.

\begin{figure}[!ht]
  \centering
  \begin{subfigure}[b]{0.4\linewidth}
    \includegraphics[width=\linewidth]{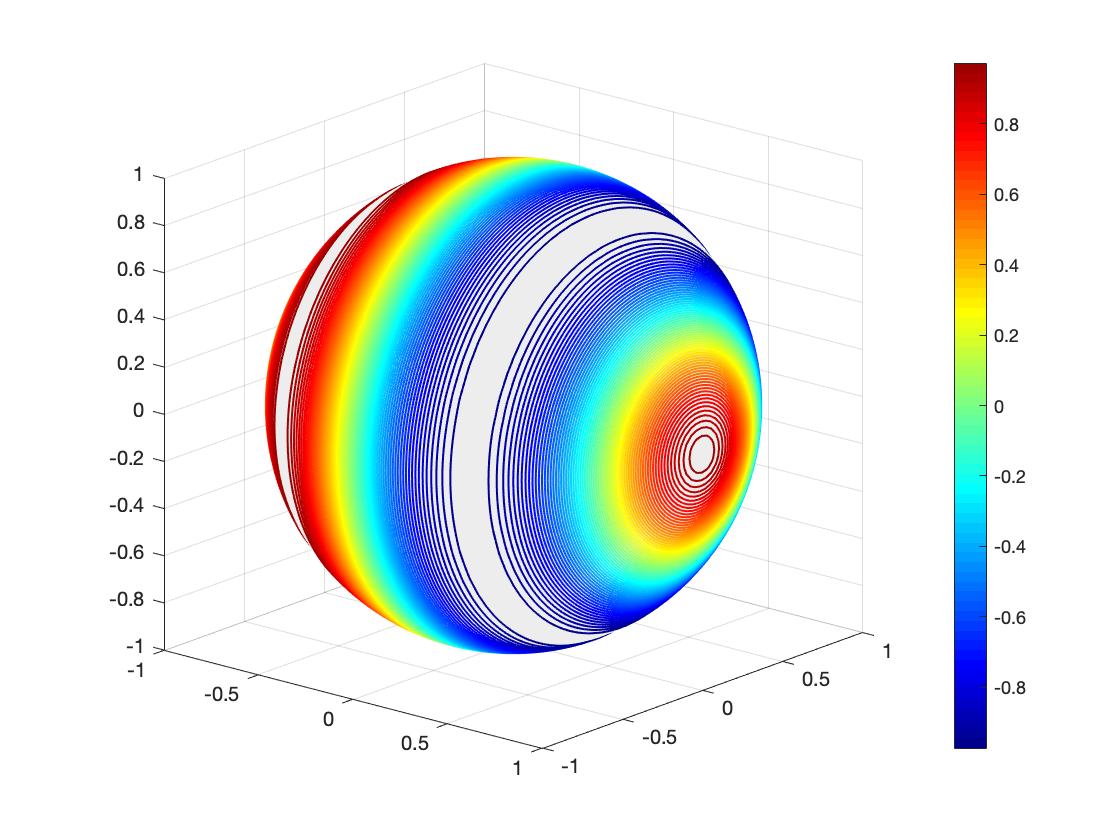}
     \caption{}
     \label{gr:case_a}
  \end{subfigure}
  \begin{subfigure}[b]{0.4\linewidth}
    \includegraphics[width=\linewidth]{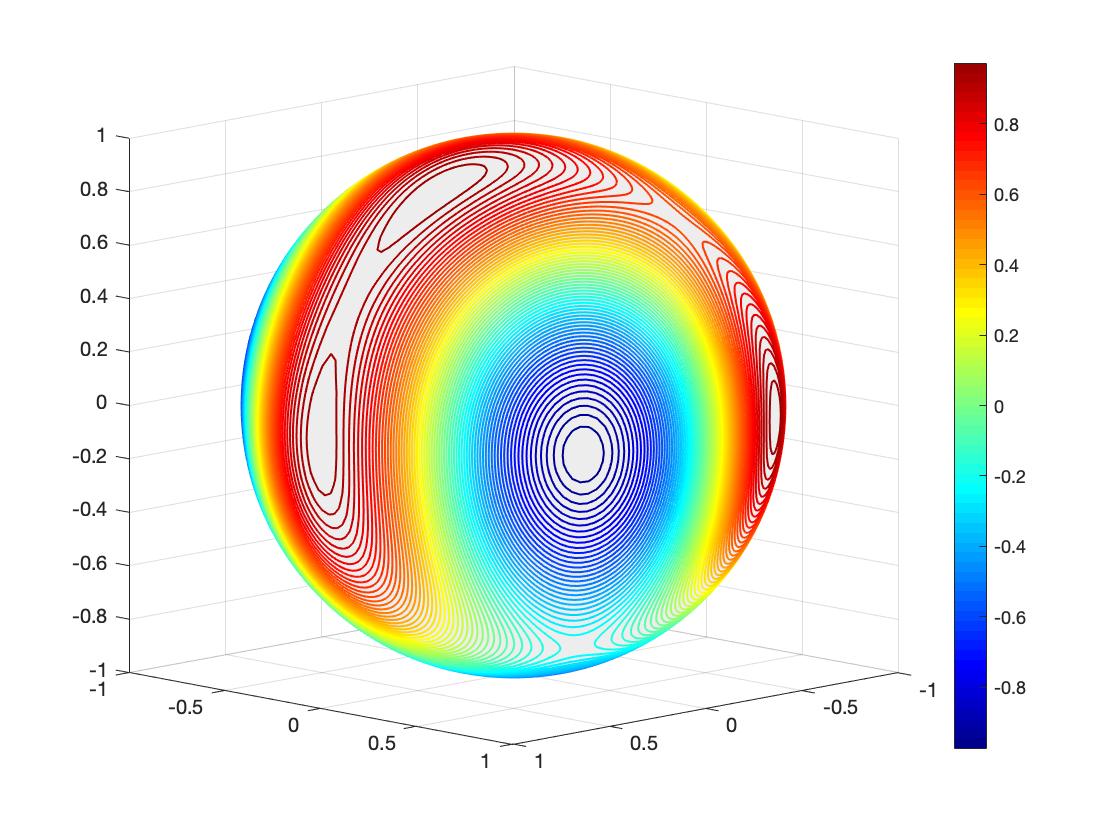}
    \caption{}
    \label{gr:case_b}
  \end{subfigure}
  \begin{subfigure}[b]{0.4\linewidth}
    \includegraphics[width=\linewidth]{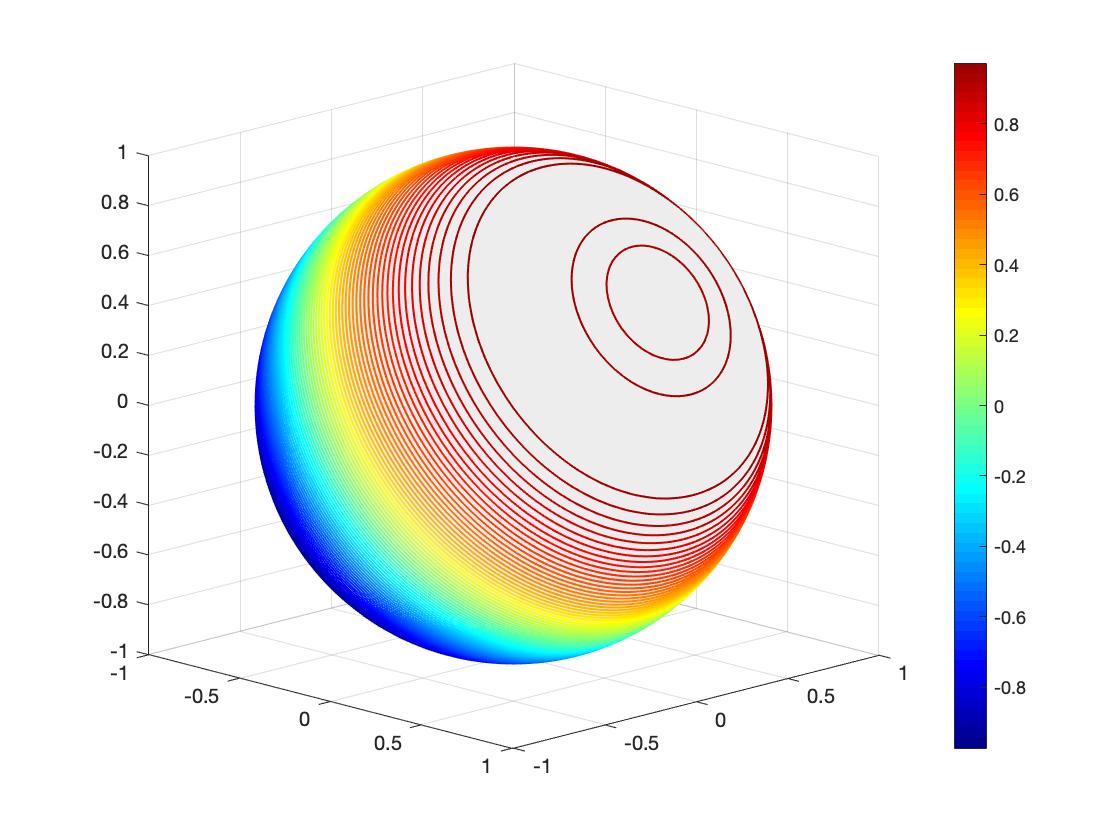}
    \caption{}
    \label{gr:case_c}
  \end{subfigure}
  \begin{subfigure}[b]{0.4\linewidth}
    \includegraphics[width=\linewidth]{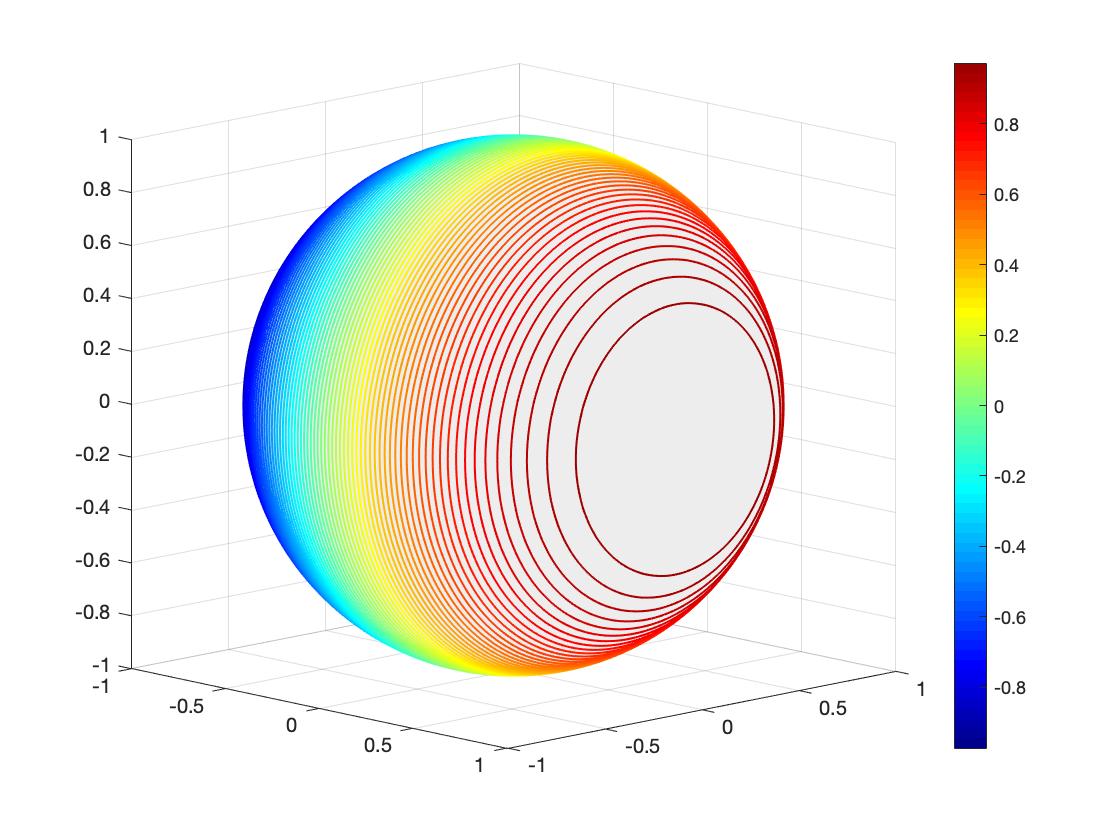}
    \caption{}
    \label{gr:case_d}
  \end{subfigure}
  \begin{subfigure}[b]{0.4\linewidth}
    \includegraphics[width=\linewidth]{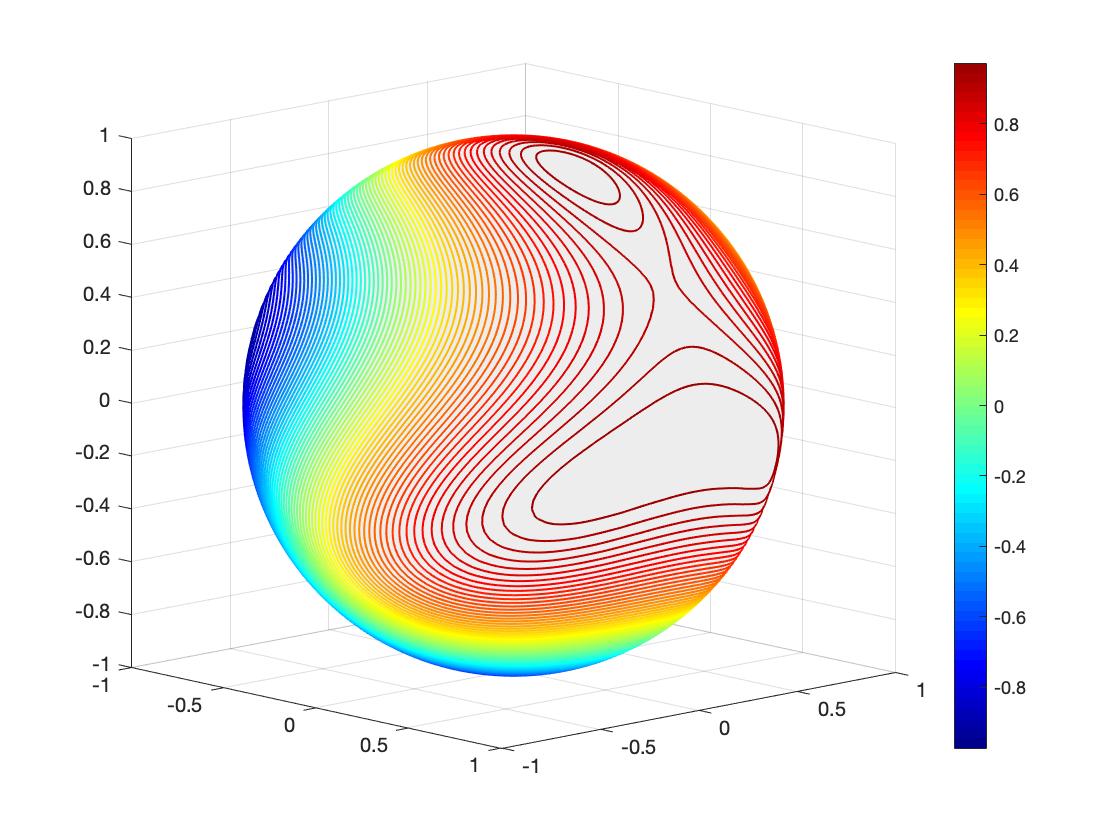}
    \caption{}
    \label{gr:case_e}
  \end{subfigure}
  \begin{subfigure}[b]{0.4\linewidth}
    \includegraphics[width=\linewidth]{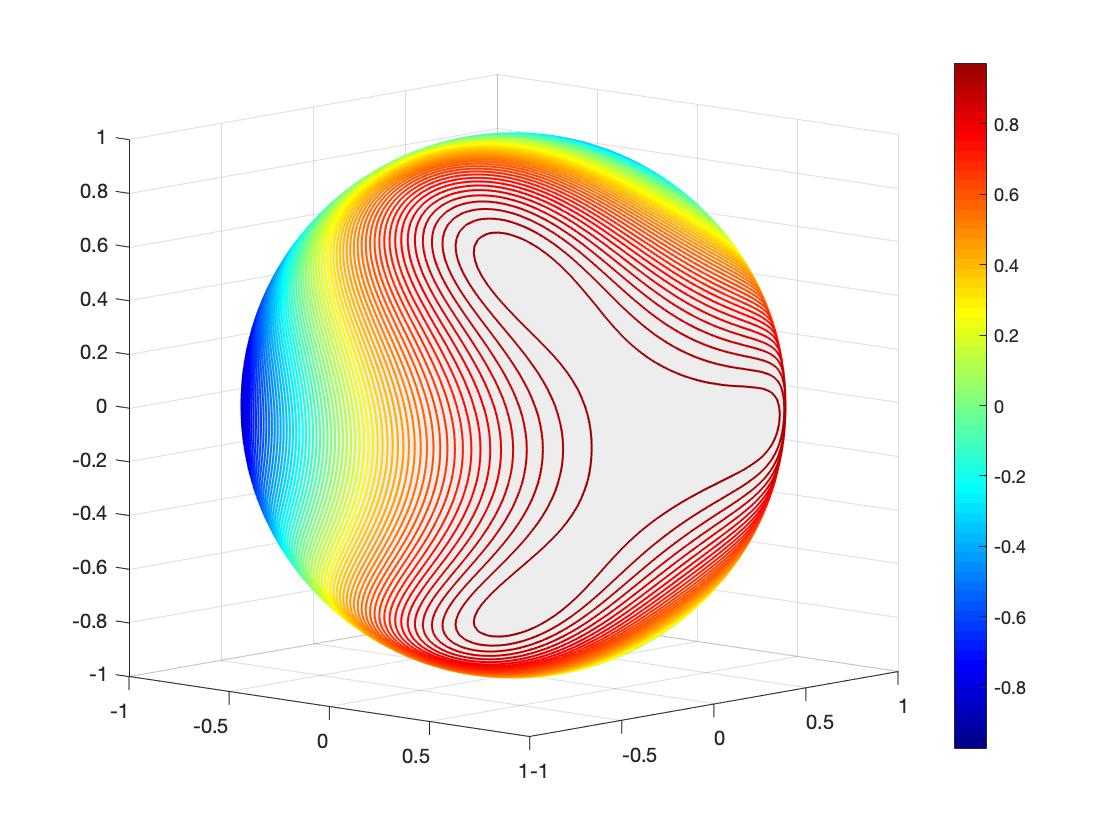}
    \caption{}
    \label{gr:case_f}
  \end{subfigure}
    \begin{subfigure}[b]{0.4\linewidth}
    \includegraphics[width=\linewidth]{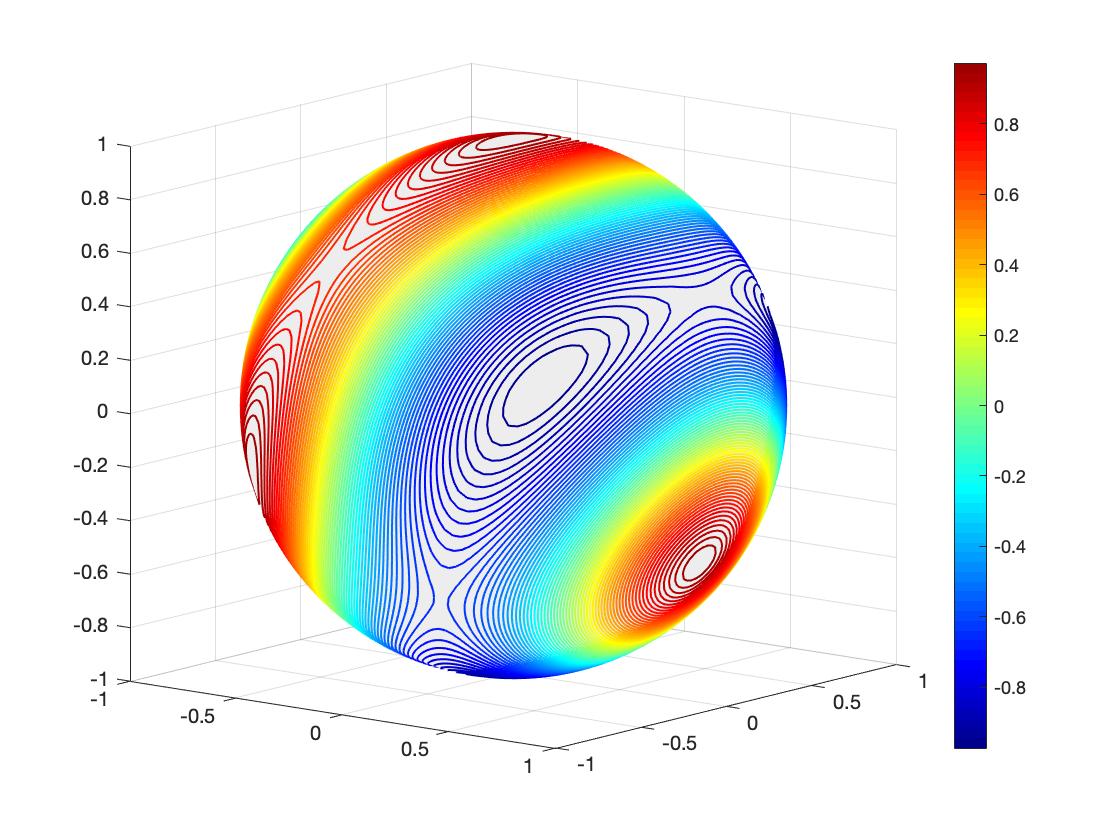}
    \caption{}
    \label{gr:case_g}
  \end{subfigure}
  \begin{subfigure}[b]{0.4\linewidth}
    \includegraphics[width=\linewidth]{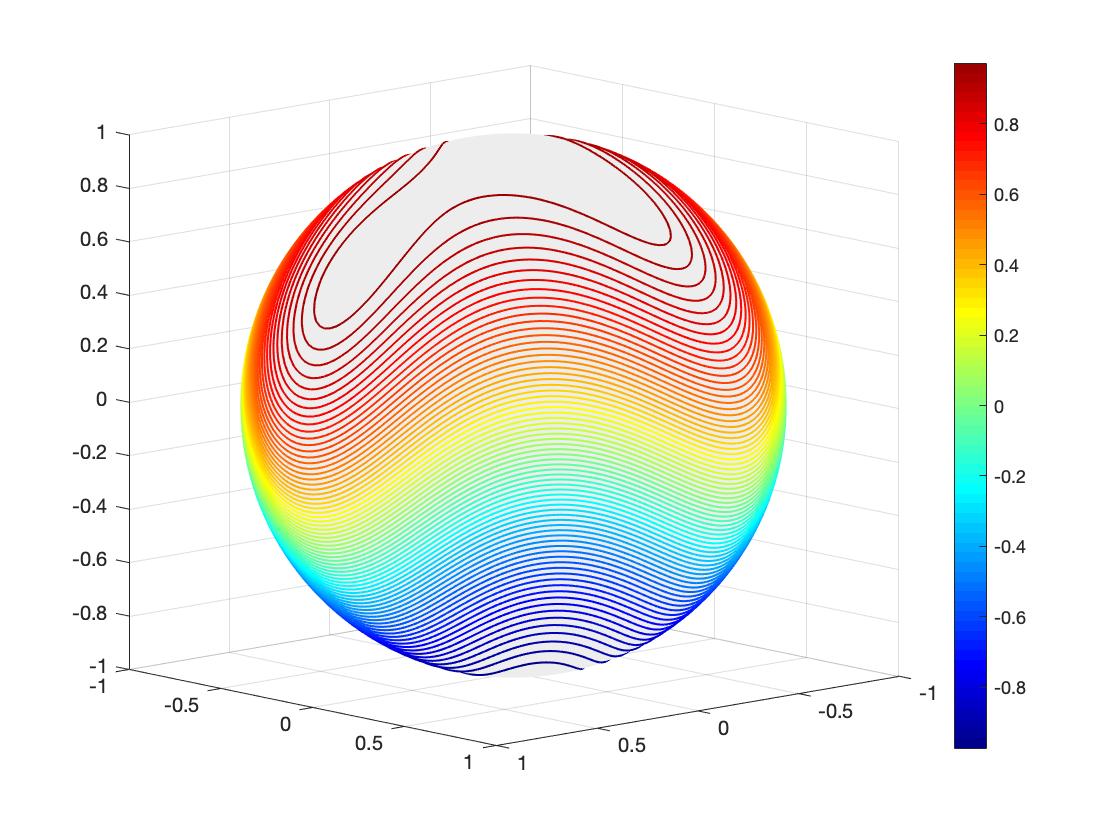}
    \caption{}
    \label{gr:case_h}
  \end{subfigure}
  \caption{Contour plots of cubic $p$ for different canonical forms.}
\end{figure}

{\theorem \label{thm:extremalS2} Let the cubic $p$ be an extremal point of the unit norm ball $B_1(S^2)$. Then either there exists an orthogonal transformation of $\mathbb R^3$ which brings $p$ to exactly one of the following canonical forms:
\begin{itemize}
\item[(a)] the zonal cubic $x_1^3 - 3x_1(x_2^2+x_3^2)$, in this case the set of global maxima of $p$ on $S^2$ consists of a non-degenerate isolated point and a circle whose points are an angle $\frac{2\pi}{3}$ apart from the isolated point (Fig. \ref{gr:case_a});
\item[(b)] the cubic $x_1^3 + 3\left(-x_2^2 + p_{102}x_3^2\right)x_1 + \left(3p_{012}x_2 + \sqrt{((1 + p_{102})^2 - 3p_{012}^2)(1 - 2p_{102})}x_3\right)x_3^2$ with $(p_{102},p_{012})$ being an arbitrary point in the triangle given by $0 \leq p_{012} \leq \sqrt{3}p_{102}$, $p_{102} < \frac12$, in this case $p$ has four isolated non-degenerate global maxima on $S^2$, three of which lie equally spaced on a great circle, and the fourth of which is an angle strictly less than $\frac{2\pi}{3}$ apart from each of the first three (Fig. \ref{gr:case_b});
\item[(c)] the zonal cubic $x_1^3 + \frac32\left(x_2^2 + 2p_{102}x_3^2\right)x_1 + \sqrt{1-2p_{102}}\left(\frac32x_2^2x_3 + (1+p_{102})x_3^3\right)$ with $p_{102} \in \left(-1,\frac12\right)$ arbitrary, in this case the global maximum on $S^2$ is achieved on some circle (Fig. \ref{gr:case_c});
\item[(d)] the zonal cubic $x_1^3 + \frac32x_1(x_2^2+x_3^2)$, in this case $e_1$ is the unique global maximum of $p$ on $S^2$, this maximum is flat (Fig. \ref{gr:case_d});
\item[(e)] the cubic $x_1^3 + 3\left(\frac12x_2^2 + p_{102}x_3^2\right)x_1 + \sqrt{1-2p_{102}}\left(-\frac32x_2^2 + (1+p_{102})x_3^2\right)x_3$, in this case $p$ has a non-degenerate global maximum at $x^+ = \frac{1}{1 - p_{102}}\left( p_{102}, 0, \sqrt{1-2p_{102}} \right)^T$ and a triply degenerate global maximum at $e_1$ (Fig. \ref{gr:case_e});
\item[(f)] the cubic \eqref{deg_nondeg_nondeg} with arbitrary $\xi \in [0,\frac{\pi}{2})$ and $p_{102} \in \left(-1,\frac12\right)$, this cubic has a doubly degenerate global maximum at $e_1$ and two non-degenerate global maxima at the points \eqref{nondeg_zeros} on $S^2$ (Fig. \ref{gr:case_f});
\end{itemize}
or $p$ has four distinct global non-degenerate maxima on $S^2$ whose Gramian is given by 
\begin{itemize}
\item[(g)] either \eqref{Gamma_central} for arbitrary $b \in \Delta_4^o$ (Fig. \ref{gr:case_g});
\item[(h)] or \eqref{Gamma_wing} for arbitrary $b \in \mathbb R_{++}^3$ (Fig. \ref{gr:case_h}).
\end{itemize} }

\begin{proof}
We consider several cases.

Suppose that $p$ has three distinct global maxima on some great circle. Then $p$ can be transformed to one of the cubics described in Corollary \ref{cor:3maxima_extremal}. The first item in this corollary yields case (a) of the theorem. The cubics described in the second item can still be transformed one into another by the coordinate transformations $x_2 \mapsto -x_2$, $x_3 \mapsto -x_3$, and rotations by $\pm\frac{2\pi}{3}$ of the plane spanned by $e_1,e_2$. These symmetries may be used to fix the sign of $p_{003}$ and to reduce the domain of the pair $(p_{102},p_{012})$ to the triangle defined in case (b).

Suppose that $p$ has a degenerate global maximum. By an orthogonal transformation of $\mathbb R^3$ $p$ can be brought to the form \eqref{p_degenerate_maximum} in Lemma \ref{lem:degenerated_maxima}, and hence satisfies the conditions in one of the 5 items of this lemma. The cubics from item 1 have already been covered in case (a) of the theorem. The transformation $x_3 \mapsto -x_3$ takes the two cubics in item 2 of Lemma \ref{lem:degenerated_maxima} to each other, and we need to list only one of them in case (c). The same reasoning yields case (e) from the two cubics in item 4. Item 3 in the lemma yields case (d). The transformations $x_2 \mapsto -x_2$, $x_3 \mapsto -x_3$ are equivalent to the substitutions $\xi \mapsto \pi - \xi$, $\xi \mapsto -\xi$ in item 5 of Lemma \ref{lem:degenerated_maxima}, reducing the set of values assumed by $\xi$ to the interval in case (f) of the theorem.

Finally, if $p$ does not have degenerate global maxima, and no three of the non-degenerate global maxima lie on a great circle, then $p$ is of the form described in Corollary \ref{cor:4global_max_necessary}, yielding cases (g),(h) of the theorem. By Lemma \ref{lem:local_to_global} the parameter $b$ can assume arbitrary values in the domains claimed in these cases.
\end{proof}

The cubics in case (b) of Theorem \ref{thm:extremalS2} have four non-degenerate global maxima on $S^2$ with Gramian
\begin{equation} \label{caseBgramian}
\begin{pmatrix} 1 & -\frac12 & -\frac12 & b_1 \\ -\frac12 & 1 & -\frac12 & b_2 \\ -\frac12 & -\frac12 & 1 & b_3 \\ b_1 & b_2 & b_3 & 1 \end{pmatrix},
\end{equation}
where $b_i \in \left( -\frac12,1 \right)$ and $b_1 + b_2 + b_3 = 0$, thus forming a triangle in the space of parameters $b = (b_1,b_2,b_3)$.

A closer look reveals that the matrices \eqref{caseBgramian}, \eqref{Gamma_central}, \eqref{Gamma_wing} and their images under conjugation with permutation matrices form a single connected 3-dimensional manifold in the space of real symmetric $4 \times 4$ matrices. The corresponding extremal cubics corresponding to cases (b),(g),(h) of Theorem \ref{thm:extremalS2} form a 6-dimensional connected manifold. The cubics corresponding to cases (a),(d),(e),(f) lie on the boundary of this manifold, while the zonal cubics corresponding to case (c) form a 3-dimensional connected manifold of extremal cubics linking the cubics corresponding to case (a) to those corresponding to case (d).

A generic extremal point of $B_1(S^2)$ can thus be of two types: either a zonal cubic with its global maximum on $S^2$ achieved on some circle, or a cubic having four non-degenerate isolated global maxima on $S^2$. The latter type further divides into two generic sub-types corresponding to cases (g),(h) in Theorem \ref{thm:extremalS2}, and the cubics corresponding to case (b) which lie between the generic sub-types. All other extremal cubics can be obtained as limits of cubics of generic types.

\bibliographystyle{plain}
\bibliography{convexity,polynomial_optim,interior_point}

\begin{thebibliography}{10}

\bibitem{AhmedStill19}
Faizan Ahmed and Georg Still.
\newblock Maximization of homogeneous polynomials over the simplex and the
  sphere: Structure, stability, and generic behavior.
\newblock {\em J. Optimiz. Theory App.}, 181:972--996, 2019.

\bibitem{BlekhermanEtAl15}
Grigoriy Blekherman, Sadik Iliman, and Martina Kubitzke.
\newblock Dimensional differences between faces of the cones of nonnegative
  polynomials and sums of squares.
\newblock {\em Int. Math. Res. Not.}, 2015:8437--8470, 2015.

\bibitem{BlekhermanParriloBook}
Grigoriy Blekherman, Pablo~A. Parrilo, and Rekha~R. Thomas, editors.
\newblock {\em Semidefinite Optimization and Convex Algebraic Geometry}.
\newblock MOS-SIAM series on Optimization. SIAM, 2013.

\bibitem{BuchheimFampaSarmiento19}
Christoph Buchheim, Marcia Fampa, and Orlando Sarmiento.
\newblock Tractable relaxations for the cubic one-spherical optimization
  problem.
\newblock In {\em Optimization of Complex Systems: Theory, Models, Algorithms
  and Applications}, volume 991 of {\em Advances in Intelligent Systems and
  Computing}, pages 267--276. Springer, 2019.

\bibitem{ChoiLam77}
Man-Duen Choi and Tsit-Yuen Lam.
\newblock Extremal positive semidefinite forms.
\newblock {\em Math. Ann.}, 231:1--18, 1977.

\bibitem{deKlerkLaurent20online}
Etienne de~Klerk and Monique Laurent.
\newblock Convergence analysis of a {L}asserre hierarchy of upper bounds for
  polynomial minimization on the sphere, 2020.
\newblock Published online in Math. Program.

\bibitem{FangFawzi21}
Kun Fang and Hamza Fawzi.
\newblock The sum-of-squares hierarchy on the sphere and applications in
  quantum information theory.
\newblock {\em Math. Program.}, 190:331--360, 2021.

\bibitem{Hilbert}
David Hilbert.
\newblock {\"U}ber die {D}arstellung definiter {F}ormen als {S}umme von
  {F}ormenquadraten.
\newblock {\em Mathematische Annalen}, 32:342--350, 1888.

\bibitem{Hildebrand21}
Roland Hildebrand.
\newblock Optimal step length for the {N}ewton method: Case of self-concordant
  functions.
\newblock https://arxiv.org/abs/2003.08650, 2021.
\newblock Accepted at Math. Methods. Oper. Res.

\bibitem{Hildebrand21b}
Roland Hildebrand.
\newblock Semi-definite representations for sets of cubics on the 2-sphere.
\newblock https://arxiv.org/abs/2103.13270, 2021.

\bibitem{KunertThesis}
Aaron Kunert.
\newblock {\em Facial structure of cones of nonnegative forms}.
\newblock PhD thesis, University Konstanz, Konstanz, 2014.

\bibitem{Naldi14}
Simone Naldi.
\newblock Nonnegative polynomials and their {C}arath\'eodory number.
\newblock {\em Discrete \& Computational Geometry}, 51:559--568, 2014.

\bibitem{NesterovSOS}
Yuri Nesterov.
\newblock Squared functional systems and optimization problems.
\newblock In Hans Frenk, Kees Roos, T\'amas Terlaky, and Shuzhong Zhang,
  editors, {\em High Performance Optimization}, chapter~17, pages 405--440.
  Kluwer Academic Press, Dordrecht, 2000.

\bibitem{NesterovMatEllipsoid}
Yuri Nesterov.
\newblock Random walk in a simplex and quadratic optimization over convex
  polytopes.
\newblock Discussion paper 2003/71, CORE, Louvain-la-Neuve, 2003.

\bibitem{Nie12}
Jiawang Nie.
\newblock Sum of squares methods for minimizing polynomial forms over spheres
  and hypersurfaces.
\newblock {\em Front. Math. China}, 7:321--346, 2012.

\bibitem{Reznick78}
Bruce Reznick.
\newblock Extremal {PSD} forms with few terms.
\newblock {\em Duke Math. J.}, 45(2):363--374, 1978.

\bibitem{Reznick00}
Bruce Reznick.
\newblock Some concrete aspects of {H}ilbert's 17th problem.
\newblock {\em Contemporary Mathematics}, 253:251--272, 2000.

\bibitem{Saunderson19}
James Saunderson.
\newblock Certifying polynomial nonnegativity via hyperbolic optimization.
\newblock {\em SIAM J. Appl. Algebra Geom.}, 3(4):661--690, 2019.

\bibitem{ManChoSo11}
Anthony Man-Cho So.
\newblock Deterministic approximation algorithms for sphere constrained
  homogeneous polynomial optimization problems.
\newblock {\em Math. Program.}, 129:357--382, 2011.

\bibitem{ZhangQiYe12}
Xinzhen Zhang, Liqun Qi, and Yinyu Ye.
\newblock The cubic spherical optimization problems.
\newblock {\em Math. Comput.}, 81:1513--1525, 2012.

\end{thebibliography}

\end{document}